# NON-ARCHIMEDEAN FLAG DOMAINS AND SEMISTABILITY I

Harm Voskuil

**Introduction.**

Let $K$ be non-archimedean local field. Then $K$ is a finite extension of $\mathbb{Q}_p$ or $\mathbb{F}_q((t))$. Let $K^\circ$ be the valuation ring of $K$. We consider a simply connected absolutely almost simple linear algebraic group $G$ defined over $K^\circ$. We always assume that $G(K)$ is not compact. For a parabolic subgroup $P \subset G$ the corresponding projective homogeneous variety is denoted by $X := G/P$. The variety $X$ cannot always be defined over $K^\circ$. For convenience we will assume in the introduction that it is defined over $K^\circ$.

Let $\mathcal{L}$ be an ample line bundle on $X$. Again we will assume in this introduction that the embedding in a projective space determined by $\mathcal{L}$ is defined over $K^\circ$. For each maximal $K$-split torus $S_K \subset G \otimes_{K^\circ} K$ we take the $S_K$-linearization of $\mathcal{L} \otimes_{K^\circ} K$ obtained by restricting the unique $G \otimes_{K^\circ} K$-linearization of $\mathcal{L} \otimes_{K^\circ} K$ to $S_K$. In this paper we will study the analytic subspaces $Y^s$ and $Y^{ss}$ of $X \otimes_{K^\circ} \mathbb{K}$. Here $\mathbb{K}$ denotes a complete non-archimedean field that contains the completion of the separable algebraic closure of $K$ and that is such that the additive valuation $v$ satisfies $v(\mathbb{K}^*) = \mathbb{R}$. We fix such a field $\mathbb{K}$. The space $Y^s$ consists of the points in $X \otimes_{K^\circ} \mathbb{K}$ that are stable for all maximal $K$-split tori $S_K \subset G \otimes_{K^\circ} K$. Similarly, $Y^{ss}$ consists of the points that are semistable for all maximal $K$-split tori.

Let $\mathcal{B}$ denote the affine building of $G(K)$. We construct a $G(K)$-equivariant map $I$, that associates to each point $x \in Y^{ss}$ a subset $I(x) \subset \mathcal{B}$. We briefly describe how the map $I$ is defined. The subset $I(x) \subset \mathcal{B}$ contains the points $z \in \mathcal{B}$ such

1991 *Mathematics Subject Classification.* 14G22, 14L24, 14M15, 20G25, 32P05.
Part of this work was supported by the Australian Research Council.

Typeset by $\mathcal{A}_{\mathcal{M}}\mathcal{S}$-TeX

that $x \in Y_z^{ss}$. Here $Y_z^{ss} \subset X \otimes_{K^\circ} \mathbb{K}$ is the open analytical subset that consists of the points $x \in X \otimes_{K^\circ} \mathbb{K}$ that are semistable in the reduction for all "relevant" maximal $K$-split tori. So for each point $z \in \mathcal{B}$ we need a reduction map. To define this reduction map we use a $\mathbb{K}^\circ$-module $V_z \subset V \otimes_{K^\circ} \mathbb{K}$. Here $V$ is the $G$-module defined over $K^\circ$ that is determined by the (very) ample line bundle $\mathcal{L}$, that is, $\mathcal{L}$ gives an embedding of $X$ into $\mathbb{P}(V)$. We will make this sketch somewhat more precise. Let us first describe, how one associates to each point $z \in \mathcal{B}$ a $\mathbb{K}^\circ$-module $V_z \subset V \otimes_{K^\circ} \mathbb{K}$.

We take a maximal $K^\circ$-split torus $S \subset G$. Then $S_K := S \otimes_{K^\circ} K \subset G \otimes_{K^\circ} K$ is a maximal $K$-split torus. The apartment $A \subset \mathcal{B}$ that corresponds to $S$ is identified with $\mathcal{X}_*(S) \otimes \mathbb{R}$. The action of $S(\mathbb{K})$ on $A$ is given by $s \cdot z = z + \nu(s)$, where $z \in A$ and $\nu(s)$ for $s \in S(\mathbb{K})$ is defined by $<\chi, \nu(s)> = -v(\chi(s))$ for all $\chi \in \mathcal{X}(S)$. Here $<,>$ is the perfect pairing between characters and one-parameter subgroups extended to a pairing between $\mathcal{X}(S) \otimes \mathbb{R}$ and $\mathcal{X}_*(S) \otimes \mathbb{R}$.

To each point $z \in A$, one associates a $\mathbb{K}^\circ$-module $V_z \subset V \otimes_{K^\circ} \mathbb{K}$. One puts $V_0 := V \otimes_{K^\circ} \mathbb{K}^\circ$ and if $z = s \cdot 0$ with $s \in S(\mathbb{K})$, then we put $V_z := s \cdot V_0$. Since $S(\mathbb{K}^\circ)$ stabilises both the point $z \in A$ and the module $V_z$, the module $V_z$ does not depend on the choice of the element $s \in S(\mathbb{K})$ with $z = s \cdot 0$. The stabiliser of $V_z$ in $G(K)$ is $P_z := s \cdot G(\mathbb{K}^\circ) \cdot s^{-1} \cap G(K)$. This is exactly the parahoric subgroup of $G(K)$ that stabilises the point $z \in A$. If $z \in \mathcal{B}$ is a point, then there exists an element $g \in G(K)$ and a point $z' \in A$ such that $g(z') = z$. Then we take $V_z := g(V_{z'})$. Since the stabiliser of $V_{z'}$ is the parahoric subgroup $P_{z'}$, the $\mathbb{K}^\circ$-module $V_z$ does not depend on the choice of the pair $g \in G(K)$ and $z' \in A$.

For every point $z \in \mathcal{B}$ we have a unique $\mathbb{K}^\circ$-module $V_z \subset V \otimes_{K^\circ} \mathbb{K}$. This gives for each point $z \in \mathcal{B}$ a unique reduction map $\psi_z : \mathbb{P}(V \otimes_{K^\circ} \mathbb{K}) = \mathbb{P}(V_z \otimes_{\mathbb{K}^\circ} \mathbb{K}) \longrightarrow \mathbb{P}(V_z \otimes_{\mathbb{K}^\circ} \overline{\mathbb{K}})$. Here $\overline{\mathbb{K}}$ denotes the residue field of $\mathbb{K}$. For each point $z \in \mathcal{B}$ one now defines an analytical subset $Y_z^{ss} \subset X \otimes_{K^\circ} \mathbb{K}$. The subset $Y_z^{ss}$ consists of the points $x \in X \otimes_{K^\circ} \mathbb{K}$ for which $\psi_z(x)$ is semistable for all tori $S \otimes_{K^\circ} \overline{\mathbb{K}}$, such that $S \subset G$ is a maximal $K^\circ$-split torus belonging to an apartment $A$ that contains the point $z$.

We show that $Y^{ss} = \bigcup_{z \in \mathcal{B}} Y_z^{ss}$. Therefore $I(x) := \{z \in \mathcal{B} \mid x \in Y_z^{ss}\}$ is a non-empty subset of the building. The subset $I(x) \subset \mathcal{B}$ is, in fact, convex. Moreover, the subset $I(x)$ is bounded if and only if $x \in Y^s$.

If the line bundle $\mathcal{L}$ is such that $Y^s = Y^{ss}$, then $I(x)$ consists of a single point for all points $x \in Y^{ss}$. In particular, this is the case when $G(K) = SL_n(K)$ and $X \otimes_{K^\circ} K = \mathbb{P}_K^{n-1}$. Then $Y^s = Y^{ss}$ is Drinfeld's symmetric space. In this particular case the map $I$ coincides with the map defined by Drinfel'd in [D] §6.

If $Y^s = Y^{ss}$, then the quotient $Y^s/\Gamma$ is proper for a discrete co-compact subgroup $\Gamma \subset G(K)$ (See [PV] and [V.1]). If $Y^s \neq Y^{ss}$, then the quotient $Y^s/\Gamma$ is still separated, but it is not proper anymore.

We use the map $I$ to construct a compactification of $Y^s$. Let $S \subset G$ be a maximal $K^\circ$-split torus and let $\chi \in \mathcal{X}(S)$ be a character. We take the $S$-linearization of $\mathcal{L}^{\otimes n}$, $n >> 0$ that one obtains by translating the $S$-linearization of $\mathcal{L}^{\otimes n}$, coming from the $G$-linearization, by the character $\chi$. For $n >> 0$ the thus obtained varieties of stable and semistable points do not depend on the particular value of $n$.

Let $P^\chi \subset G(K)$ be the parabolic subgroup that is generated by the Borel subgroups $B \subset G(K)$ that contain $S_K := S \otimes_{K^\circ} K$ and that are such that the character $\chi$ is contained in the positive Weyl chamber belonging to $B$. All the maximal $K$-split tori of $G(K)$ that are contained in $P^\chi$ are conjugated in $P^\chi$. Hence the $S_K$-linearization of $\mathcal{L}^{\otimes n} \otimes_{K^\circ} K$, $n >> 0$ induces an $S'_K$-linearization of $\mathcal{L}^{\otimes n} \otimes_{K^\circ} K$ for every maximal $K$-split torus $S'_K \subset P^\chi$. Let $Y^s_\chi$ (resp. $Y^{ss}_\chi$) be the analytic subvariety of $X \otimes_{K^\circ} \mathbb{K}$ that consists of the points that are stable (resp. semistable) for all maximal $K$-split tori $S'_K \subset P^\chi$ for our linearization of $\mathcal{L}^{\otimes n} \otimes_{K^\circ} \mathbb{K}$, $n >> 0$. Then $Y^s \subseteq Y^s_\chi \subseteq Y^{ss}_\chi \subseteq Y^{ss}$.

There exist characters $\chi \in \mathcal{X}(S)$ such that $Y^s_\chi = Y^{ss}_\chi$. Such a character $\chi$ gives a compactification of $Y^s$ in the following sense. There exists a formal scheme $\mathcal{Y}_\chi$ over $spf(\mathbb{K}^\circ)$ with generic fibre $Y^s_\chi$ and whose closed fibre consists of proper components that correspond 1-1 to the vertices of the building $\mathcal{B}$ of $G(K)$. The group $G(K)$ does not act on $Y^s_\chi$ nor on $\mathcal{Y}_\chi$. However, the parabolic subgroup $P^\chi \subset G(K)$ does act on both $Y^s_\chi$ and $\mathcal{Y}_\chi$.

Our terminology differs a little bit from [PV] and [V.1]. We call $Y^s$ a non-archimedean flag domain. In [V.1] one calls $Y^s$ a flag domain if and only if $Y^s/\Gamma$ is proper. In [PV] such spaces were called symmetric spaces. However it seems better to preserve the name symmetric spaces for those spaces $Y^s$ with automorphic functions for $\Gamma$ on it and such that $X = G/P$ with $P$ maximal parabolic. The spaces $Y^{ss}$ are called period domains in [Ra] and [RZ].

From the point of view of rigid analytic uniformization the space $Y^{ss}$ itself is not very useful, since a discrete co-compact subgroup $\Gamma \subset G(K)$ does not act discretely on $Y^{ss}$ if $Y^{ss} \neq Y^s$. However from a moduli perspective it is quite interesting. (See [Ra], [To] and [RZ]).

The first three sections of this paper are a revision of the first three sections of [V.2].

## §1. Root systems, ample line bundles and semistability

**1.1.** In this section we collect some preliminaries from the litterature. We describe the relative and absolute root system of $G$. The connection between dominant weights and ample line bundles is described. Using this one can find the field over which a non-archimedean flag domain can be defined. We also recall some criteria with which one can determine whether a point is (semi-)stable. Our references for the group theory quoted here are [Hu], [S] and [T.1].

**1.2. The absolute root system.** Let $G$ be a simply connected absolutely almost simple linear algebraic group defined over $K^\circ$. Let $K_s$ denote the separable closure of $K$ and $K_s^\circ$ its ring of integers.

Let $S \subseteq T \subset G$ be such that $S \otimes_{K^\circ} K$ is a maximal $K$-split torus and $T$ a maximal torus defined over $K^\circ$. The torus $T$ and the group $G$ both split over $K_s$. Let $\mathcal{X}(T)$ be the character group of $T$. Let $\Phi$ be the *absolute root system* of $G$. We choose a *simple basis* $\Delta$ of $\Phi$.

The *Weyl group* of $\Phi$ is denoted by $W$. The group $W$ acts on $\mathcal{X}(T) \otimes \mathbb{R}$. After choosing a $W$-invariant inner product on $\mathcal{X}(T) \otimes \mathbb{R}$, the Weyl group $W$ is generated by the reflections in the hyperplanes orthogonal to the simple roots $\alpha \in \Delta$. Let $N_T$ and $Z_T$ be the *normaliser* and *centraliser* of the torus $T$ in $G$. Then $W \cong N_T/Z_T(K_s)$.

Let $\alpha_j, j = 1, \ldots, \ell$ be the simple roots in $\Delta$. Here $\ell := \dim_{K_s}(T \otimes_{K^\circ} K_s)$ is the *absolute rank* of $G$. Let $w_j \in W$ denote the reflection in the hyperplane orthogonal to $\alpha_j$.

The simple basis $\Delta$ determines a *Borel subgroup* $B \subset G$. Let $\Phi^+$ denote the set of *positive roots* for $\Delta$. To each character $\alpha \in \Phi$ belongs a $T$-stable additive subgroup $U_\alpha \subset G$ on which $T$ acts with character $\alpha$. One has $B = <T, U_\alpha \mid \alpha \in \Phi^+>$.

Let $J \subseteq \{1, \ldots, \ell\}$ be a non-empty subset. The subgroup $W_J \subseteq W$ is the group generated by the reflections $w_j, j \notin J$. The *parabolic subgroup* $P_J \supseteq B$ is the group $P_J := BW_JB$. Any parabolic subgroup of $G$ is conjugated to exactly one of the groups $P_J$.

**1.3. The relative root system.** We fix an ordering on $\mathcal{X}(S)$ and choose a compatible ordering on $\mathcal{X}(T)$. Let $\Delta_0$ denote the set of simple roots of $G$ with respect to $T$ that vanish on $S$. The *relative root system* $\Phi_K$ of $G$ consists of the roots of $\mathcal{X}(S)$. The *relative Weyl group* $W_K$ acts on $\Phi_K$. One has $W_K \cong N_S/Z_S(K)$. Here $N_S$ and $Z_S$ are the normaliser and centraliser of $S$ in $G$.

The Galois group $H := Gal(K_s/K)$ acts on $\mathcal{X}(T)$, since $T$ is split over $K_s$. For any element $h \in H$ the image $h(\Delta)$ of $\Delta$ is again a simple basis of $\Phi$. Hence there exists a unique element $w \in W$ such that $wh(\Delta) = \Delta$. We put $h^* := w \circ h$ for this element $w \in W$. We call this action of $H$ on $\mathcal{X}(T)$ the *twisted action*. Let $H^* := \{h^* | h \in H\}$. Then $H^*$ is a finite group acting on $\mathcal{X}(T)$ that preserves the simple basis $\Delta$.

**1.4. Ample line bundles.** Every ample line bundle $\mathcal{L}$ on $X := G/P$ is very ample. So $\mathcal{L}$ determines an embedding $X \hookrightarrow \mathbb{P}(V)$. Here $\mathbb{P}(V)$ denotes the projectivisation of a $G$-module $V$. The module $V$ is irreducible if $Char(K) = 0$, whereas if $Char(K) > 0$ this might not be the case. (See [J]). The $G$-module $V$ is uniquely determined by its highest weight $\lambda$. We will denote this $G$-module by $V_\lambda$.

Let $(-,-)$ be a $W$-invariant inner product on $\mathcal{X}(T) \otimes \mathbb{R}$. The *fundamental weights* $\omega_i$ of our root system $\Phi$ are determined by $2(\omega_i, \alpha_j)/(\alpha_j, \alpha_j) = \delta_{ij}$. The ample line bundles $\mathcal{L}$ on $X = G/P_J$ correspond to the modules $V_\lambda$ with highest weight $\lambda = \sum_{j \in J} n_j \omega_j$ with $n_j > 0$ for all $j \in J$.

Let $\Lambda_r$ denote the *root lattice*, i.e. the sublattice of $\mathcal{X}(T)$ generated by the roots $\alpha \in \Phi$. In [T.2] theorems 3.3 and 7.2 the irreducible representations of $G$ that are defined over $K$ are determined. From this one gets the linear representations of $G$. Furthermore, these representations are defined over $K^\circ$.

For each $\lambda \in \Lambda_r$ one has a representation $\rho_\lambda$ of $G$ into $GL(\oplus_{\sigma \in H^*} V_{\sigma(\lambda)})$ that is defined over $K^\circ$. (If the weight $\lambda$ is not contained in $\Lambda_r$ one might need a skew field $\mathcal{D}$ defined over $K$ to define $\rho_\lambda$). For our purposes the representations $\rho_\lambda, \lambda \in \Lambda_r$ suffice.

Note that $\rho_\lambda$ is irreducible if $Char(K) = 0$, whereas if $Char(K) = p > 0$ this might not be the case.

**1.5. Fields of definition.** For each weight $\lambda = \sum_{j \in J} n_j \omega_j \in \Lambda_r$ we have a representation of $G$ defined over $K^\circ$. Let $v_\lambda \in \oplus_{\sigma \in H^*} V_{\sigma(\lambda)}$ be a vector contained in $V_\lambda \oplus <0>$, whose component in $V_\lambda$ is a highest weight vector. The image $X_\lambda$ of the orbit $G \cdot v_\lambda$ in $\mathbb{P}(\oplus_{\sigma \in H^*} V_{\sigma(\lambda)})$ is isomorphic to $G/P_J$. We have a variety $X^\dagger \subset \mathbb{P}(\oplus_{\sigma \in H^*} V_{\sigma(\lambda)})$ defined over $K^\circ$, whose connected components are $X_{\sigma(\lambda)}$, $\sigma \in H^*$. All the connected components of $X^\dagger$ are isomorphic. The very amply line bundle $\mathcal{L}^\dagger$ associated with this embedding restricts on each connected component $X_{\sigma(\lambda)}$ to the line bundle $\mathcal{L}$ belonging to the weight $\sigma(\lambda)$.

For each $X = G/P_J$ defined over $K_s^\circ$ and ample line bundle $\mathcal{L}$ on $X$ corresponding to some weight $\lambda \in \Lambda_r$, one has a variety $X^\dagger$ and an ample line bundle $\mathcal{L}^\dagger$, both defined over $K^\circ$, such that one connected component of $X^\dagger$ is isomorphic to $X$ and $\mathcal{L}^\dagger|_X \cong \mathcal{L}$.

Let $L$ be a splitting field for $G$ and let $L^\circ$ be the ring of integers of $L$. Then the connected components of $X^\dagger \otimes_{K^\circ} L^\circ$ are defined over $L^\circ$. Moreover, also the restriction of $\mathcal{L}^\dagger \otimes_{K^\circ} L^\circ$ to a connected component of $X^\dagger \otimes_{K^\circ} L^\circ$ is defined over $L^\circ$. Any field that has these properties will be called a *field of definition* of the pair $(X, \mathcal{L})$. One can find a minimal field of definition of the pair $(X, \mathcal{L})$.

**1.6. (Semi)-stable points.** Let $\mathcal{L}$ be an ample line bundle on $X$ corresponding to a weight $\lambda$, then $\mathcal{L}^{\otimes n}$ corresponds to the weight $n \cdot \lambda$. For every weight $\lambda$ one can find an integer $n > 0$ such that $n \cdot \lambda \in \Lambda_r$. Moreover, $\mathcal{L}$ and $\mathcal{L}^{\otimes n}$ determine the same sets of (semi-)stable points. So the assumption that $\lambda \in \Lambda_r$ poses no real restriction. We will from now on assume that both $X$ and $\mathcal{L}$ are defined over $L^\circ$, where $L$ is some field of definition of the pair $(X, \mathcal{L})$.

Let $X^s(T, \mathcal{L})$ and $X^{ss}(T, \mathcal{L})$ denote the varieties consisting of the stable and the semistable points, respectively, in $X$ for the action of $T \otimes_{K^\circ} L^\circ$ w.r.t. $\mathcal{L}$. The $T \otimes_{K^\circ} L^\circ$ linearization is chosen to be the restriction of the unique $G \otimes_{K^\circ} L^\circ$ linearization of $\mathcal{L}$ (or $\mathcal{L}^\dagger$). These sets can be determined using the criteria given in [MFK] Ch. 2 (See also [DH] 1.1.5).

Let $\mathcal{L}$ give an embedding of $X$ into $\mathbb{P}(V)$ for some $G \otimes_{K^\circ} L^\circ$-module $V := V_\lambda$. We have a decomposition $V = \oplus V^\beta, \beta \in \mathcal{X}(T)$ into eigenspaces $V^\beta$ on which $T \otimes_{K^\circ} L^\circ$ acts with character $\beta$. Each $V^\beta$ has a free basis $e_1^\beta, \cdots, e_{m(\beta)}^\beta$, defined over the ring of integers of some finite separable extension $L' \supseteq L$. These bases determine coordinates $x_{\beta,i}$ for every $x \in X \subset \mathbb{P}(V)$. We may assume that $x$ is normalised, i.e. $|x_{\beta,i}| \leq 1, \forall \beta, i$ and $\max |x_{\beta,i}| = 1$. Then one defines the polyhedron $\mu(x) \subset \mathcal{X}(T) \otimes \mathbb{R}$ as the convex hull of $\{\beta \mid \exists i \; |x_{\beta,i}| = 1\}$. We also define a polyhedron $\mu_K(x) \subset \mathcal{X}(T) \otimes \mathbb{R}$ as being the convex hull of $\{\beta \mid \exists i \; x_{\beta,i} \neq 0\}$. One has for $x \in X$:

$$\bar{x} \in X^s(T, \mathcal{L}) \otimes_{L^\circ} \bar{L} \iff 0 \in \operatorname{int} \mu(x)$$

$$\bar{x} \in X^{ss}(T, \mathcal{L}) \otimes_{L^\circ} \bar{L} \iff 0 \in \mu(x)$$

$$x \in X^s(T, \mathcal{L}) \otimes_{L^\circ} L \iff 0 \in \operatorname{int} \mu_K(x)$$

$$x \in X^{ss}(T, \mathcal{L}) \otimes_{L^\circ} L \iff 0 \in \mu_K(x)$$

Here $\bar{L}$ denotes the residue field of $L$ and $\bar{x}$ denotes the reduction of $x$. Part a) of the proposition below is proved in [GS] and [FH]. Part b) was proved in [PV] theorem 1.1:

**1.7. Proposition.** *Let $\mathcal{L}$ be an ample line bundle on $X = G/P$ corresponding to the weight $\lambda$. Let $T \subset G$ be a maximal torus, defined over $K^\circ$. Then one has:*

a) *For any point $x \in X$ the vertices of $\mu(x)$ form a subset of the set $\{w(\lambda) | w \in W\}$. The edges of $\mu(x)$ are parallel to the roots $\alpha \in \Phi$.*

b) *$X^s(T, \mathcal{L}) = X^{ss}(T, \mathcal{L})$ if and only if $\lambda$ is not contained in a hyperplane (through 0) spanned by roots.*

**1.8.** Let $S$ be a maximal $K^\circ$-split torus and let $T \supseteq S$ be a maximal torus defined over $K^\circ$. Let $r : \mathcal{X}(T) \longrightarrow \mathcal{X}(S)$ denote the restriction map. The varieties of (semi-) stable points $X^s(S, \mathcal{L})$ and $X^{ss}(S, \mathcal{L})$ for the torus $S \otimes_{K^\circ} L^\circ$ can be described using the polyhedron $r(\mu(x)) \subset \mathcal{X}(S) \otimes \mathbb{R}$. Below we recall some results and their proof from [V.1].

**1.9. Proposition.** *Let $\mathcal{L}$ be an ample line bundle on $X = G/P$, corresponding to the weight $\lambda \in \Lambda_r$ and let $S \subseteq T \subset G$ as before. Then one has:*

a) *For any point $x \in X$ the vertices of $r(\mu(x))$ form a subset of $\{r(w(\lambda)) | w \in W\}$. The edges of $r(\mu(x))$ are parallel to roots $\alpha \in \Phi_K$.*

b) $X^s(S, \mathcal{L}) = X^{ss}(S, \mathcal{L})$ *if for all $w \in W$ the weight $r(w(\lambda))$ is not contained in a hyperplane (through 0) spanned by roots $\alpha \in \Phi_K$.*

c) $X^s(S, \mathcal{L}) \neq X^{ss}(S, \mathcal{L})$ *if $\lambda$ is contained in a hyperplane $V \subset \mathcal{X}(T) \otimes \mathbb{R}$ spanned by roots $\alpha \in \Phi$ and moreover there exists an element $w \in W$ such that $r(w(V))$ is contained in a hyperplane (spanned by roots $\beta \in \Phi_K$).*

*Proof.* Since $r$ maps $\Phi$ into $\Phi_K \cup \{0\}$, part (a) of the proposition is clear. Furthermore (b) follows from (a).

To prove (c) one constructs for $\lambda$ in a hyperplane $V$ spanned by roots $\alpha \in \Phi$ a point $x \in X$ such that $\mu(x) = V \cap \text{conv}(\{w(\lambda)|w \in W\})$. This has been done in [PV] 1.4. Then one chooses an element $w \in W$ such that $r(w(V))$ is contained in a hyperplane in $\mathcal{X}(S) \otimes \mathbb{R}$. Then $0 \in r(w(\mu(x))) = r(\mu(w(x))) \subset r(w(V))$ and the point $w(x)$ is in $X^{ss}(S, \mathcal{L}) - X^s(S, \mathcal{L})$. This proves (c).

**1.10.** In the theorem below we describe the projective homogeneous varieties on which there exist ample line bundles, such that the notions of semistability and stability coincide. For a proof based on a case by case study of the root systems and proposition 1.9 above, we refer to [V.1]. As usual, we assume that $G(K)$ is non-compact.

**1.11 Theorem.** *There exists an ample line bundle $\mathcal{L}$ on $X = G/P$ such that $X^{ss}(S, \mathcal{L}) = X^s(S, \mathcal{L})$ if and only if one of the following holds:*

(1) $P = B$ *and $G$ is any group.*

(2) $P = P_J$, $J = \{1, \ldots, \ell - 1\}$ *and $G$ is a non-split group with absolute root system $\Phi$ of type $C_\ell$. Here $\alpha_\ell$ is the unique long root in the simple basis $\Delta$.*

(3) $P = P_J$ *and $G = SL_{s+1}(\mathcal{D})$ with g.c.d.$(j \in J, s+1) = 1$. Here $\mathcal{D}$ is a skew field (of dimension $d^2$) defined over $K$ (we allow $\mathcal{D} = K$) and the $P_{\{j\}}$ are such that $G/P_{\{j\}} \cong Gr(j, (s+1) \cdot d)$.*

§2. EQUIVARIANT MAPS TO THE APARTMENT

**2.1.** We fix a group $G$, a projective homogeneous variety $X = G/P$ and an ample line bundle $\mathcal{L}$ on $X$ corresponding to some weight $\lambda \in \Lambda_r$. We fix a field of definition $L$ for the pair $(X, \mathcal{L})$ such that $X$ and $\mathcal{L}$ are both defined over $L^\circ$. Let $S \subset G$ be a maximal $K^\circ$-split torus and let $A$ be the apartment corresponding to $S$.

By $\mathbb{K}$ we denote a complete non-archimedean field that contains the completion of the separable closure of $K$ and that is such that $v(\mathbb{K}^*) = \mathbb{R}$. Here $v$ denotes the additive valuation of $\mathbb{K}$. We fix such a field $\mathbb{K}$. Using such a large field $\mathbb{K}$, has the advantage that the torus $S(\mathbb{K})$ acts transitively on the apartment $A$. In this section we construct an $S(\mathbb{K})$-equivariant map $I_A : X^{ss}(S, \mathcal{L}) \otimes_{L^\circ} \mathbb{K} \longrightarrow \{$ convex subsets of $A$ $\}$. The map $I_A$ will be studied in some detail. We also briefly discuss the possibility of extending the map $I_A$ to all of $X \otimes_{L^\circ} \mathbb{K}$.

First we recall the construction of the apartment $A$ belonging to $S$ and describe the action of $S(\mathbb{K})$ on $A$ (See [La] Ch. 1 §1 and [T.3] 1.1-1.2).

**2.2. The apartment.** Let $\mathcal{X}_*(S)$ be the space of one-parameter subgroups (1-ps), i.e. the dual of the character space $\mathcal{X}(S)$. The perfect pairing between characters and 1-ps will be denoted by $<-,->$. We extend it to a pairing between $\mathcal{X}(S) \otimes \mathbb{R}$ and $\mathcal{X}_*(S) \otimes \mathbb{R}$.

One defines a map $\nu : S(\mathbb{K}) \longrightarrow \mathcal{X}_*(S) \otimes \mathbb{R}$ by $<\chi, \nu(s)> = -v(\chi(s))$ for all $\chi \in \mathcal{X}(S)$. Here $v$ is the additive valuation of $\mathbb{K}$, normalized such that $v(\pi) = 1$ for a generator $\pi$ of the maximal ideal of $K^\circ$. Hence $v(K^*) = \mathbb{Z}$. We regard $\mathcal{X}_*(S) \otimes \mathbb{R}$ as an affine space. The torus $S$ acts on $\mathcal{X}_*(S) \otimes \mathbb{R}$ by translation, as follows: $s \cdot z = z + \nu(s)$ for $z \in \mathcal{X}_*(S) \otimes \mathbb{R}$. This gives a well-defined action of $S(\mathbb{K})$ on $\mathcal{X}_*(S) \otimes \mathbb{R}$. Since $v(\mathbb{K}^*) = \mathbb{R}$, the torus $S(\mathbb{K})$ acts transitively on the points in the affine space $\mathcal{X}_*(S) \otimes \mathbb{R}$.

To $G$ and $S$ belongs an affine root system $\Phi_{aff}$. The affine roots are functions $\alpha + n : \mathcal{X}_*(S) \otimes \mathbb{R} \longrightarrow \mathbb{R}$ with $\alpha \in \Phi_K$ and $n \in \Gamma_\alpha$. Here $(\alpha + n)(z) = \alpha(z) + n$ for $z \in \mathcal{X}_*(S) \otimes \mathbb{R}$, where $\alpha(z) := <\alpha, z>$. Furthermore $\Gamma_\alpha$ is the discrete group $\mathbb{Z}$ or $\frac{1}{2}\mathbb{Z}$, depending on the particular affine root system $\Phi_{aff}$ and the root $\alpha$.

The affine root system gives a simplicial decomposition of $\mathcal{X}_*(S) \otimes \mathbb{R}$. The maximal simplices or *chambers* are the closures of the connected components of $\{z \in \mathcal{X}_*(S) \otimes \mathbb{R} \mid \forall (\beta \in \Phi_{aff}) \quad \beta(z) \neq 0\}$. The affine space $\mathcal{X}_*(S) \otimes \mathbb{R}$ endowed with this simplicial decomposition is the *apartment $A$* for the torus $S \otimes_{K^\circ} K$.

**2.3. Analytification of the varieties of (semi-) stable points.** To each algebraic variety corresponds a rigid analytic variety that has the same set of closed points. (See [BGR] or [FP]). We denote the analytifications of $X^s(S, \mathcal{L}) \otimes_{L^\circ} \mathbb{K}$ and $X^{ss}(S, \mathcal{L}) \otimes_{L^\circ} \mathbb{K}$ by $Y_A^s$ and $Y_A^{ss}$, respectively. For notational convenience we have dropped the line bundle from the notation (though these spaces do depend on $\mathcal{L}$). The subscript corresponds to the apartment $A$ that belongs to $S(K)$. From now on we will also regard the projective algebraic variety $X \otimes_{L^\circ} \mathbb{K}$ as a rigid analytic variety.

To each point $z \in A$, we want to associate an analytic subspace $Y^s_{z,A}$ (resp. $Y^{ss}_{z,A}$) of $Y^s_A$ (resp. $Y^{ss}_A$). One defines $Y^s_{0,A}$ (resp. $Y^{ss}_{0,A}$) as being the completion of $X^s(S, \mathcal{L}) \otimes_{L^\circ} \mathbb{K}^\circ$ (resp. $X^{ss}(S, \mathcal{L}) \otimes_{L^\circ} \mathbb{K}^\circ$) along the closed fiber. Here $\mathbb{K}^\circ$ is the ring of integers of $\mathbb{K}$. One has by construction:

$$Y^s_{0,A}(\mathbb{K}) = X^s(S, \mathcal{L})(\mathbb{K}^\circ)$$
$$Y^{ss}_{0,A}(\mathbb{K}) = X^{ss}(S, \mathcal{L})(\mathbb{K}^\circ) \,.$$

Let $z \in A$ be a point and let $s \in S(\mathbb{K})$ be such that $z = s \cdot 0$. Then we define:

$$Y^s_{z,A} := s \cdot Y^s_{0,A} = \{s \cdot x \mid x \in Y^s_{0,A}\}$$
$$Y^{ss}_{z,A} := s \cdot Y^{ss}_{0,A} = \{s \cdot x \mid x \in Y^{ss}_{0,A}\} \,.$$

Since $S(\mathbb{K}^\circ)$ is the stabiliser in $S(\mathbb{K})$ of both the point $0 \in A$ and of the subspaces $Y^s_{0,A}, Y^{ss}_{0,A} \subset Y^{ss}_A$, the definitions of the spaces $Y^s_{z,A}$ and $Y^{ss}_{z,A}$ do not depend on the choice of the element $s \in S(\mathbb{K})$ such that $z = s \cdot 0$. The spaces $Y^s_{z,A}$ and $Y^{ss}_{z,A}$ are both quasi-compact, i.e. they are the union of finitely many affinoid subspaces.

Another way to define the spaces $Y^s_{z,A}$ and $Y^{ss}_{z,A}$ is as follows. To the point $0 \in A$ one associates the $\mathbb{K}^\circ$-module $V_0 := V \otimes_{L^\circ} \mathbb{K}^\circ \subset V \otimes_{L^\circ} \mathbb{K}$. For $z \in A$ one takes $V_z := s \cdot V_0$, where $s \in S(\mathbb{K})$ is such that $s \cdot 0 = z$. Since $S(\mathbb{K}^\circ)$ is the stabiliser of both the point $0 \in A$ and the module $V_0$, the module $V_z$ is well-defined and unique. For each point $z \in A$ we now have a reduction map $\psi_z : \mathbb{P}(V \otimes_{L^\circ} \mathbb{K}) = \mathbb{P}(V_z \otimes_{\mathbb{K}^\circ} \mathbb{K}) \longrightarrow \mathbb{P}(V_z \otimes_{\mathbb{K}^\circ} \overline{\mathbb{K}})$. Then $Y^{ss}_{z,A}$ (resp. $Y^s_{z,A}$) consists of the points $x \in X \otimes_{L^\circ} \mathbb{K}$ such that $\psi_z(x)$ is semistable (resp. stable) for the action of $S \otimes_{K^\circ} \overline{\mathbb{K}}$.

**2.4. Proposition.**

i) $Y^{ss}_{z,A} \subset Y^{ss}_A$ and $Y^s_{z,A} \subset Y^s_A$.

ii) $Y^{ss}_A = \bigcup_{z \in A} Y^{ss}_{z,A}$.

iii) $Y^s_A = \bigcup_{z \in A} Y^s_{z,A}$ if and only if $Y^s_A = Y^{ss}_A$.

*Proof.* The first statement is obvious. So let us consider statement (ii).

Since $X^{ss}(S, \mathcal{L})/(S \otimes_{K^\circ} L^\circ)$ is a projective algebraic variety, the map $\psi^{ss} : S(\mathbb{K}) \times Y_{0,A}^{ss} \longrightarrow Y_A^{ss}$, defined by $\psi^{ss}(s,x) = s \cdot x$, is surjective. From this statement (ii) follows.

Statement (iii) follows from the fact that the map $\psi^s : S(\mathbb{K}) \times Y_{0,A}^s \longrightarrow Y_A^s$, given by $(s,x) \longrightarrow s \cdot x$, is surjective if and only if $Y_A^s = Y_A^{ss}$.

**2.5. Definition.** For $x \in Y_A^{ss}$ we define the *interval of S-semistability* $I_A(x)$ as follows:

$$I_A(x) := \{z \in A \mid x \in Y_{z,A}^{ss}\}$$
$$= \{s^{-1} \cdot 0 \in A \mid s \cdot x \in Y_{0,A}^{ss}\}$$

Note that $I_A(s \cdot x) = s \cdot I_A(x)$ for all $s \in S(\mathbb{K})$ and $x \in Y_A^{ss}$. Since $Y_A^{ss}$ is the union of the subspaces $Y_{z,A}^{ss}$, $z \in A$, the subset $I_A(x) \subseteq A$ is non-empty for all $x \in Y_A^{ss}$.

**2.6. Proposition.** *Let $x \in Y_A^{ss}$ be a point. Then:*

a) $I_A(x) \subseteq A$ *is convex.*

b) $I_A(x)$ *is bounded if and only if $x \in Y_A^s$.*

c) $I_A(x) = \{z\}$ *if and only if $x \in Y_{z,A}^s$.*

*Proof.* Since $I_A(s \cdot x) = s \cdot I_A(x)$ for $s \in S(\mathbb{K})$, we may assume that $0 \in I_A(x)$, i.e. that $x \in Y_{0,A}^{ss}$. We assume that $x \in Y_{0,A}^{ss}$ is normalised. Let the coordinates $x_{\beta,i}$ be as in 1.6 and let $r : \mathcal{X}(T) \longrightarrow \mathcal{X}(S)$ be the restriction map.

One has : $s \cdot x \in Y_{0,A}^{ss}$ if and only if $|r(\beta)(s) \cdot x_{\beta,i}| \leq 1 \quad \forall \beta, i$. Indeed, this follows from the fact that $s$ preserves (the absolute value of) the $S$-invariants. So if $s \in S(\mathbb{K})$ is such that $s \cdot x \in Y_{0,A}^{ss}$, then $s$ satisfies:

$$-v(r(\beta)(s)) \leq v(x_{\beta,i}) \quad \forall \quad \beta, i \quad s.t. \quad x_{\beta,i} \neq 0 .$$
$$\iff -v(r(\beta)(s)) \leq \min \{v(x_{\beta,i}) \mid x_{\beta,i} \neq 0\} \quad \forall \beta \quad s.t. \quad \exists i \text{ with } x_{\beta,i} \neq 0 .$$

Let $V_K(x) \subset \mathcal{X}(T)$ consist of the characters $\beta \in \mathcal{X}(T)$ such that some coordinate $x_{\beta,i} \neq 0$. For $\beta \in V_K(x)$ we put $n_\beta := min \{v(x_{\beta,i}) \mid x_{\beta,i} \neq 0\}$. Then $z = s^{-1} \cdot 0 \in I_A(x)$ if and only if the point $z$ satisfies for all $\beta \in V_K(x)$:

$$<r(\beta), z> = <r(\beta), s^{-1} \cdot 0> = <r(\beta), -\nu(s)> = v(r(\beta)(s)) \geq -n_\beta.$$

Therefore $I_A(x) \subseteq A$ is the intersection of a finite number of affine halfspaces $H(r(\beta), n_\beta) := \{z \in A \mid <r(\beta), z> \geq -n_\beta\}$. In particular, $I_A(x)$ is convex. This proves a).

Let $\chi \in \mathcal{X}(S) \otimes \mathbb{R}$. If $\chi = \sum f_\beta \cdot r(\beta)$ with $f_\beta \in \mathbb{R}_{\geq 0}$ and $\beta \in V_K(x)$, then $<\chi, z> = \sum f_\beta \cdot <r(\beta), z> \geq -\sum f_\beta \cdot n_\beta$ for all $z \in I_A(x)$. If $\chi$ is not contained in $F_x^+ := \{\sum f_\beta \cdot r(\beta) \mid f_\beta \in \mathbb{R}_{\geq 0}, \beta \in V_K(x)\}$, then $<\chi, z>$ is not bounded from below for $z \in I_A(x)$. Hence $I_A(x)$ is bounded if and only if $F_x^+ = \mathcal{X}(S) \otimes \mathbb{R}$. Furthermore, $F_x^+ = \mathcal{X}(S) \otimes \mathbb{R}$ if and only if 0 is contained in the convex hull of $r(V_K(x))$. The convex hull of $r(V_K(x))$ is $r(\mu_K(x))$. Therefore $I_A(x)$ is bounded if and only if $x \in Y_A^s$. This proves b).

The proof of c) resembles that of b) and is left to the reader.

**2.7. Proposition.** *The interval of semistability $I_A(x)$ for $x \in Y_A^{ss}$ is bounded by hyperplanes parallel to walls. In particular, there exist $n_\alpha(x) \in \mathbb{R} \cup \{\infty\}, \alpha \in \Phi_K$ such that:*

$$I_A(x) = \{z \in A \mid \forall (\alpha \in \Phi_K) \quad \alpha(z) \leq n_\alpha(x)\} .$$

*Proof.* Without loss of generality, we may assume that $x \in Y_{0,A}^{ss}$. Furthermore, we may assume that 0 is not contained in a face of the polyhedron $I_A(x)$. We normalise $x$ such that $\max \{|x_{\beta,i}| \mid \beta \in \mathcal{X}(T)\} = 1$. Let $V_1(x, A) \subset \mathcal{X}(S) \otimes \mathbb{R}$ be the linear subspace spanned by the weights $r(\beta), \beta \in \mathcal{X}(T)$ such that for some $i$ one has $|x_{\beta,i}| = 1$. From proposition 1.9 it follows that either $V_1(x, A) = <0>$ or $V_1(x, A)$ is spanned by roots $\alpha \in \Phi_K$.

If $V_1(x, A) = \mathcal{X}(S) \otimes \mathbb{R}$ we are done. Indeed, then $x \in Y_{0,A}^s$ and $I_A(x) = \{0\}$. Hence $n_\alpha(x) = 0$ for all $\alpha \in \Phi_K$.

So let us assume that $V_1(x, A) \neq \mathcal{X}(S) \otimes \mathbb{R}$. Let $S_1 \subset S$ be the subtorus on which all characters $\chi \in V_1(x, A) \cap \mathcal{X}(S)$ vanish. Clearly, if $s \in S(\mathbb{K})$ is such that $s \cdot x \in Y_{0,A}^{ss}$, then $s \in S_1(\mathbb{K})$. Hence $I_A(z) \subseteq \{z \in A \mid \alpha(z) = 0 \quad \forall (\alpha \in V_1(x, A) \cap \Phi_K)\}$.

Let $s \in S_1(\mathbb{K})$ be such that $s^{-1} \cdot 0$ is in the interior of a codimension one face of $I_A(x)$. Let $F \subseteq \mathcal{X}(S) \otimes \mathbb{R}$ be the linear subspace spanned by the $r(\beta)$ such that some $|x_{\beta,i}(s \cdot x)| = 1$. Then by proposition 1.9 a) the space $F$ is generated by roots $\alpha \in \Phi_K$. Moreover, $\dim (F) = \dim (V_1(x, A)) + 1$. Hence one can find a root $\alpha \in F \cap \Phi$ that together with $V_1(x, A)$ spans $F$. The face of $I_A(x)$ containing $s^{-1} \cdot 0$ is the intersection of $I_A(x)$ with a hyperplane defined by $\alpha(z) = \alpha(s^{-1} \cdot 0)$.

Furthermore, $I_A(x)$ is either contained in the half space $\{z \in A|\quad \alpha(z) \leq \alpha(s^{-1}\cdot 0)\}$ or in the halfspace $\{z \in A|\quad \alpha(z) \geq \alpha(s^{-1}\cdot 0)\}$. Since we can do this for any codimension one face of $I_A(x)$, the proposition follows.

**2.8. Remark.** Note that $I_A(x) = A$ if and only if $x \in Y_A^{ss}$ is a fixed point for the action of $S(\mathbb{K})$. Such a fixed point exists if and only if $r(w(\lambda)) = 0$ for some $w \in W$. In particular, this never occurs if $G$ is split.

If none of the characters $r(w(\lambda)), w \in W$ vanishes on $S$, then $I_A(x) \subseteq A$ has always some codimension $> 0$.

**2.9. Compactification of the apartment.** In order to extend the map $I_A$ to the unstable points of the analytic variety $X \otimes_{L^\circ} \mathbb{K}$, we compactify the apartment $A$. This is done by adding the *spherical apartment* $A_\infty$ to $A$ (See [B] Ch. VI §9 or [Ro] Ch. 9 §3). We briefly recall the construction.

A *ray* is a halfline starting at some point $z \in A$. Two rays are called *parallel* if and only if they are translations of eachother. The points of the spherical apartment $A_\infty$ are the equivalence classes that consist of parallel rays. The compactification $\overline{A}$ of the apartment $A$ consists of the union of $A$ with $A_\infty$ with the usual topology. One calls $A_\infty$ the *spherical apartment at infinity*.

As a representative of an equivalence class of rays, one can take the ray starting at the point $0 \in A$. Then one can identify $A_\infty$ with $(\mathcal{X}(S)_* \otimes \mathbb{R} - \{0\})/\sim$. Here the equivalence relation is given by $\delta_1 \sim \delta_2$ if and only if there exists an $n \in \mathbb{R}_{>0}$ such that $\delta_1 = n \cdot \delta_2$.

The spherical apartment $A_\infty$ has a simplicial structure. The simplices of $A_\infty$ are given by the equivalence classes of rays such that the representative of the class that starts at $0 \in A$ intersects a fixed simplex of $A$ that contains $0$ in a point different from $0$.

**2.10. Definition.** Let $\delta_1$, $\delta_2$ be two points in $A_\infty$. Then the points $\delta_1$ and $\delta_2$ are called *antipodal* if the corresponding rays are opposite. This means that the rays corresponding to the points $\delta_1$ and $\delta_2$ that start in the same point $z \in A$ form a line together. A subset $F \subset A_\infty$ is called *semi-convex* if it contains the line joining any pair of points $\delta_1$ and $\delta_2$ in $F$, provided that the points are not antipodal.

The *semi-convex hull* of a subset $F \subset A_\infty$ is the smallest semi-convex subset of $A_\infty$ that contains $F$.

For a point $x \in X \otimes_{L^\circ} \mathbb{K}$ we put $F_A(x) := \{\delta \in \mathcal{X}_*(S) \mid <\chi, \delta> \leq 0 \text{ for all } \chi \in r(\mu_K(x))\}$. By $\widetilde{F}_A(x)$ we denote the image of $F_A(x)$ in $A_\infty \cong (\mathcal{X}(S)_* \otimes \mathbb{R} - \{0\})/\mathbb{R}_{>0}$. Then we denote the semi-convex hull of $\widetilde{F}_A(x)$ by $\Lambda_A(x)$. Note that $\Lambda_A(x)$ depends very much on the line bundle $\mathcal{L}$.

Now we can define a map $\overline{I_A} : X \otimes_{L^\circ} \mathbb{K} \to \{\text{subsets of } \overline{A}\}$ as follows: $\overline{I_A}(x) := I_A(x) \cup \Lambda_A(x)$. Here we take $I_A(x) = \emptyset$ for $x \in X \otimes_{L^\circ} \mathbb{K} - Y_A^{ss}$.

**2.11. Proposition.** *Let $x \in X \otimes_{L^\circ} \mathbb{K}$ be a point.*

   i) $\Lambda_A(x) = \emptyset$ *if and only if* $x \in Y_A^s$.
   
   ii) *If* $x \in Y_A^{ss}$ *then the points of* $\Lambda_A(x)$ *correspond to the rays that are contained in* $I_A(x)$.

*Proof.* This is clear from the definitions.

**2.12. Remark.** The proposition above shows that $\overline{I_A}(x)$ is the closure of $I_A(x)$ in $\overline{A}$ if $x \in Y_A^{ss}$. However, for the points $x \in X \otimes_{L^\circ} \mathbb{K}$ that are not contained in $Y_A^{ss}$ the map $\overline{I_A}$ is not that good.

Indeed, there exists ample line bundles $\mathcal{L}$ and $\mathcal{L}'$ on $X$ such that $X^{ss}(S, \mathcal{L}) = X^{ss}(S, \mathcal{L}')$. One easily sees that both line bundles determine the same interval of $S$-semistability. However, if the weights $\lambda$ and $\lambda'$ that correspond to $\mathcal{L}$ and $\mathcal{L}'$ are such that $n \cdot \lambda \neq m \cdot \lambda'$ for all $n, m \in \mathbb{Z}_{>0}$, then these line bundles give rise to different maps $\overline{I_A}$. Indeed, then the maps differ on the unstable points. So the map $\overline{I_A}$ is not uniquely determined by the map $I_A$.

## §3 Equivariant maps to the building

**3.1.** Let $S \subset G$ be a maximal $K^\circ$-split torus and let $A$ be the apartment that corresponds to $S_K := S \otimes_{K^\circ} K$. Let $Y^s := \bigcap_{g \in G(K)} g(Y_A^s)$ and $Y^{ss} := \bigcap_{g \in G(K)} g(Y_A^{ss})$.

We will construct a $G(K)$-equivariant map $I$ from $Y^{ss}$ to the set of convex subsets of $\mathcal{B}$. Here $\mathcal{B}$ is the Bruhat-Tits building of $G(K)$.

To each point $z \in \mathcal{B}$ we associate some analytic subspaces $Y_z^s, Y_z^{ss} \subset X \otimes_{L^\circ} \mathbb{K}$. One proves that $Y_z^s \subset Y^s$ and $Y_z^{ss} \subset Y^{ss}$. We show that $Y^{ss} = \bigcup_{z \in \mathcal{B}} Y_z^{ss}$. Then one defines $I(x)$ as being $I(x) := \{z \in \mathcal{B} \mid x \in Y_z^{ss}\}$. We show that $I(x) \subset \mathcal{B}$ is convex. Furthermore, $I(x)$ is bounded if and only if $x \in Y^s$. We also briefly discuss how to extend the map $I$ to all of $X \otimes_{L°} \mathbb{K}$. Then the image is a subset of $\overline{\mathcal{B}}$. Here $\overline{\mathcal{B}}$ is the compactification of the affine building $\mathcal{B}$ obtained by adding as a boundary the spherical building of $G(K)$.

**3.2. The Building.** All maximal $K$-split tori $S'_K \subset G \otimes_{K°} K$ are conjugated. To each maximal $K$-split torus $gS_K g^{-1}, g \in G(K)$, belongs an apartment $g(A)$. Each simplex $\sigma \in A$ is stabilised by a parahoric subgroup $P_\sigma \subset G(K)$. The simplex $g(\sigma), g \in G(K)$, is stabilised by $P_{g(\sigma)} := gP_\sigma g^{-1}$. A point $z \in A$ is stabilised by the parahoric subgroup $P_z \subset G(K)$. Here $P_z = P_\sigma$, where $\sigma \in A$ is the smallest simplex that contains the point $z$.

On the points in the set of apartments $g(A)$, $g \in G(K)$ we define an equivalence relation $\sim$. Let $z_i \in g_i(A)$, $i = 1, 2$. Then $z_1 \sim z_2$ if and only if there exists an element $g \in P_{z_1}$ such that $z_2 = g(z_1)$ and $g_2(A) = g(g_1(A))$. Note that $P_{z_1} = P_{z_2}$ if $z_1 \sim z_2$. The *affine building* $\mathcal{B}$ is defined as $\mathcal{B} := \bigcup_{g \in G(K)} g(A)/\sim$.

**3.3. Definition.** Let $z' \in A'$. We can find an element $g \in G(K)$ such that $A' = g(A)$. We put $z := g^{-1}(z') \in A$. Then $Y^s_{z',A'} := g(Y^s_{z,A})$ and $Y^{ss}_{z',A'} := g(Y^{ss}_{z,A})$.

We associate the following analytical subspaces of $X \otimes_{L°} \mathbb{K}$ to a point $z \in \mathcal{B}$:

$$Y^s_z := \bigcap_{A' \ni z} Y^s_{z,A'} = \bigcap_{g \in P_z} g(Y^s_{z,A'}) \,.$$

$$Y^{ss}_z := \bigcap_{A' \ni z} Y^{ss}_{z,A'} = \bigcap_{g \in P_z} g(Y^{ss}_{z,A'}) \,.$$

Here $P_z \subset G(K)$ denotes the stabiliser of $z \in \mathcal{B}$.

The spaces $Y^s_z$ and $Y^{ss}_z$ are defined using an infinite number of intersections. However, it follows from proposition 3.4 below, that only a finite number of intersections are needed. Therefore the spaces $Y^s_z$ and $Y^{ss}_z$ are quasi-compact, that is, they are the union of a finite number of affinoids. In fact, if $z \in A$, then $Y^s_z \subseteq Y^s_{z,A}$ and $Y^{ss}_z \subseteq Y^{ss}_{z,A}$ are open analytical subspaces.

As in 2.3 one can use $\mathbb{K}°$-submodules of $V \otimes_{L°} \mathbb{K}$ to define the spaces $Y^s_z$ and $Y^{ss}_z$. Let us fix a maximal $\mathbb{K}°$-split torus $S \subset G$ and let $A \subset \mathcal{B}$ be the apartment that belongs to $S$. We assume that we have for each point $z \in A$ a $\mathbb{K}°$-module $V_z \subset V \otimes_{L°} \mathbb{K}$ defined as in 2.3. The stabiliser in $G(K)$ of the $\mathbb{K}°$-module $V_z$ is $s \cdot G(\mathbb{K}°) \cdot s^{-1} \cap G(K)$. Here $s \in S(\mathbb{K})$ is such that $z = s \cdot 0$. One verifies that the stabiliser equals $P_z$, the parahoric subgroup of $G(K)$ that stabilises the point $z \in A \subset \mathcal{B}$.

For any point $z \in \mathcal{B}$ there exists an element $g \in G(K)$ such that for some point $z' \in A$ one has $z = g(z')$. Then one puts $V_z := g(V_{z'})$. Since the stabiliser in $G(K)$ of $V_z$ equals the parahoric subgroup $P_z \subset G(K)$, the $\mathbb{K}^\circ$-module $V_z$ does not depend on the pair $g \in G(K)$ and $z' \in A$ such that $z = g(z')$. Therefore we have for each point $z \in \mathcal{B}$ a unique $\mathbb{K}^\circ$-module $V_z \subset V \otimes_{L^\circ} \mathbb{K}$.

For each point $z \in \mathcal{B}$, one has a reduction map $\psi_z : \mathbb{P}(V \otimes_{L^\circ} \mathbb{K}) = \mathbb{P}(V_z \otimes_{\mathbb{K}^\circ} \mathbb{K}) \longrightarrow \mathbb{P}(V_z \otimes_{\mathbb{K}^\circ} \overline{\mathbb{K}})$. Then the space $Y_z^{ss}$ (resp. $Y_z^s$) consists of the points $x \in X \otimes_{L^\circ} \mathbb{K}$, such that $\psi_z(x)$ is semistable (resp. stable) for all tori $S \otimes_{K^\circ} \overline{\mathbb{K}}$, such that $S \subset G$ is a maximal $K^\circ$-split torus corresponding to an apartment $A \subset \mathcal{B}$ that contains the point $z$.

**3.4. Proposition.** *Let $z \in \mathcal{B}$ be a point and let $A, A' \subset \mathcal{B}$ be apartments that contain the point $z$. If for some $\varepsilon > 0$ the intersection $A \cap A'$ contains the ball with center $z$ and radius $\varepsilon$, then $Y_{z,A}^s = Y_{z,A'}^s$ and $Y_{z,A}^{ss} = Y_{z,A'}^{ss}$.*

*Proof.* In the statement of the proposition, we have tacitly assumed, that we have a notion of distance on the apartment. This distance is defined by choosing some $W_K$-invariant inner product on $A$.

Let us assume that a ball with center $z$ and radius $\varepsilon > 0$ is contained in the intersection $A \cap A'$. Then for all 1-ps $\delta$ of $S$ one has that $\delta(t) \cdot z \in A \cap A'$, for $|t|$ sufficiently close (depending on $\delta$) to 1. If $g \in P_z$ is such that $g(A) = A'$, then $\delta(t) \cdot g \cdot \delta(t)^{-1} \in s \cdot G(\mathbb{K}^\circ) \cdot s^{-1}$ for all 1-ps $\delta$ of $S$ and for $|t|$ sufficiently close to 1. Here $s \in S(\mathbb{K})$ is such that $s \cdot 0 = z$.

Therefore the reduction $\overline{g}$ of $g$ in $s \cdot G \cdot s^{-1} \otimes_{K^\circ} \overline{\mathbb{K}^\circ}$ is contained in the centraliser of $S \otimes_{K^\circ} \overline{\mathbb{K}^\circ}$. Note that $s \cdot G \cdot s^{-1}$ is defined over $\mathbb{K}^\circ$ in such a way that its reduction is isomorphic with $G \otimes_{K^\circ} \overline{\mathbb{K}^\circ}$. Since the centraliser of the torus preserves both the set of the stable points and the set of the semistable points, it follows that $g$ preserves both $Y_{z,A}^s$ and $Y_{z,A}^{ss}$. This proves the proposition. 

**3.5. Definition.** Let $f_1, \ldots, f_m$ be an $L^\circ$-basis of the $L^\circ$-module $\Gamma(X, \mathcal{L}^{\otimes d})^{S \otimes_{K^\circ} L^\circ}$, consisting of the $S \otimes_{K^\circ} L^\circ$-invariants. Here $d$ is chosen in such a way that the homogeneous $S \otimes_{K^\circ} L^\circ$-invariants of degree $d$ generate $\oplus_{n>0} \Gamma(X, \mathcal{L}^{\otimes dn})^{S \otimes_{K^\circ} L^\circ}$ as a $L^\circ$-algebra.

Let $A_1 = g_1(A)$ and $A_2 = g_2(A)$ be two different apartments, where $g_1, g_2 \in$

$G(K)$. As in [PV] section 3.6, we define a function $r_{A_1,A_2} : Y^{ss}_{A_2} \longrightarrow \mathbb{R}$ as follows:

$$r_{A_1,A_2}(x) = r_{g_1(A),g_2(A)}(x) := \max_{1 \leq i \leq m}\{|g_1^* f_i(x)|\} / \max_{1 \leq i \leq m}\{|g_2^* f_i(x)|\} .$$

This function has the following obvious properties:

a) $r_{A_1,A_2}(x) = r_{A_2,A_1}(x)^{-1}$, if $x \in Y^{ss}_{A_1} \cap Y^{ss}_{A_2}$.

b) $r_{g(A_1),g(A_2)}(g(x)) = r_{A_1,A_2}(x)$ for $g \in G(K)$ and $x \in Y^{ss}_{A_2}$.

Below we will study the map $r_{A_1,A_2}$ in more detail. We follow [PV] closely here.

**3.6. Lemma.** *Let $x \in Y^{ss}_{z,A}$ and assume that $z \in A \cap A'$. Then $r_{A',A}(x) \leq 1$.*

*Proof.* Since $z$ is contained in $A \cap A'$, there exists an element $g \in P_z$ such that $A' = g(A)$. Here $P_z \subset G(K)$ is the stabiliser of $z \in A$. The lemma now follows from the fact that $P_z = s \cdot G(\mathbb{K}^\circ) \cdot s^{-1} \cap G(K)$, where $s \in S(\mathbb{K})$ is such that $s \cdot 0 = z$. (See [PV] 3.6 property c).

**3.7. Definition.** For $x \in X \otimes_{L^\circ} \mathbb{K}$ we define the function $\widetilde{r} : X \otimes_{L^\circ} \mathbb{K} \longrightarrow \mathbb{R}$ by:

$$\widetilde{r}(x) = \begin{cases} 0 & \text{if } x \notin Y^{ss}_A \\ \inf\{r_{gA,A}(x) \mid g \in G(K)\} & \text{if } x \in Y^{ss}_A \end{cases}$$

**3.8. Proposition.** *Let $x \in Y^{ss}_A$ be a point. Then there exists an apartment $A' \subset \mathcal{B}$ such that $r_{A',A}(x) = \widetilde{r}(x)$*

*Proof.* See [PV] 3.6, proof of property d).

**3.9. Proposition.** *Let $z \in A' \subset \mathcal{B}$ and let $x \in Y^{ss}_{z,A'}$. Then $x \in Y^{ss}_z$ if and only if $r_{A',A}(x) = \widetilde{r}(x) > 0$.*

*Proof.* We first prove that if $\widetilde{r}(x) = 0$, then $x \notin Y^{ss}_z$. If $\widetilde{r}(x) = 0$ then by proposition 3.8 we can find an apartment $A_1 \subset \mathcal{B}$ such that $r_{A_1,A}(x) = 0$. Hence $x \notin Y^{ss}_{A_1}$. Now we can find $z_1 \in A_1$ and an apartment $A_2$ such that $x \in Y^{ss}_{z_1,A_2}$. Since $A_1 = hA_2$ for some $h \in P_{z_1}$ we can construct $\widetilde{h} \in P_z$ such that $r_{\widetilde{h}A_2,A_2}(x) < \varepsilon$. In particular we may assume that $0 < r_{\widetilde{h}A_2,A}(x) < r_{A',A}(x)$. Now $x \in Y^{ss}_{z_2,\widetilde{h}A_2}$ for some $z_2 \in \widetilde{h}A_2$.

Let $A''$ be such that $z_2, z \in A''$. Using lemma 3.6, one finds that $r_{A'',A}(x) < r_{A',A}(x)$. In particular, it follows that $x \notin Y^{ss}_{z,A''}$. Hence $x \notin Y^{ss}_z$.

Now let us assume that $\widetilde{r}(x) > 0$. Let $A_1 \subset \mathcal{B}$ be an apartment such that $r_{A_1,A}(x) = \widetilde{r}(x)$. Suppose $x \in Y^{ss}_{z_1,A_1}$. Let $A_2$ be an apartment containing both $z$ and $z_1$. Using lemma 3.6, one finds $r_{A_2,A}(x) = \widetilde{r}(x) > 0$. Again applying 3.6 gives us: $x \in Y^{ss}_z \iff r_{A',A}(x) = \widetilde{r}(x) > 0$.

**3.10. Corollary.** *Let $z \in \mathcal{B}$. Then $Y_z^{ss} \subset Y^{ss}$.*

*Proof.* If $x \in Y_z^{ss}$ then $\widetilde{r}(x) > 0$ by proposition 3.9. Hence for all $A' \subset \mathcal{B}$ one has $x \in Y_{A'}^{ss}$. Thus $x \in Y^{ss}$.

**3.11. Theorem.** $Y^{ss} = \bigcup_{z \in \mathcal{B}} Y_z^{ss}$.

*Proof.* In the corollary above we have already shown that $Y_z^{ss} \subset Y^{ss}$. Therefore it is sufficient to prove that $Y^{ss} \subseteq \bigcup_{z \in \mathcal{B}} Y_z^{ss}$.

Let $x \in Y^{ss}$ be a point. Then there exists an apartment $A' \subset \mathcal{B}$ such that $r_{A',A}(x) = \widetilde{r}(x)$. Hence $x \in Y_{z,A'}^{ss}$ for some $z \in A'$. Using proposition 3.9, one concludes that $x \in Y_z^{ss}$. This proves the theorem.

**3.12. Proposition.** *Let $x \in X \otimes_{L^\circ} \mathbb{K}$. Then $x \in Y^{ss}$ if and only if $\widetilde{r}(x) > 0$.*

*Proof.* Clear from proposition 3.9 and theorem 3.11.

**3.13. Definition.** For $x \in Y^{ss}$ one defines the *interval of $G(K)$-semistability $I(x)$* as follows:
$$I(x) := \{z \in \mathcal{B} \mid x \in Y_z^{ss}\} \ .$$

Since $Y_z^{ss} = \bigcap_{A' \ni z} Y_{z,A'}^{ss}$ one has:
$$I(x) = \{z \in \mathcal{B} \mid \forall (A' \ni z) \quad z \in I_{A'}(x)\} \ .$$

Note that $I_{gA}(x) = g(I_A(g^{-1}(x)))$ for $g \in G(K)$. Clearly, one has by construction that $I(g(x)) = g(I(x))$ for $g \in G(K)$. It follows from theorem 3.11 above, that $I(x)$ is non-empty for any point $x \in Y^{ss}$.

**3.14. Theorem.** *Let $x \in Y^{ss}$ be a point. Then the following statements hold:*

i) $I(x) = \bigcup \{I_{A'}(x) \mid r_{A',A}(x) = \widetilde{r}(x)\}$.
ii) *If $A' \subset \mathcal{B}$ is an apartment such that $A' \cap I(x) \neq \emptyset$, then $I_{A'}(x) = A' \cap I(x)$.*
iii) *The interval of $G(K)$-semistability $I(x)$ is convex.*

*Proof.* The first statement follows immediately from the definition of $I(x)$ and propositions 3.9 and 3.12.

Let us now consider the second statement. If $z \in A' \cap I(x)$, then $x \in Y_z^{ss}$. In particular, $z \in I_{A'}(x)$ and therefore $A' \cap I(x) \subseteq I_{A'}(x)$. Moreover, $r_{A',A}(x) = \widetilde{r}(x) > 0$ and therefore $I_{A'}(x) \subseteq I(x)$. This proves the second statement.

To prove the third statement, we consider two points $z_1, z_2 \in I(x)$. If $A' \subset \mathcal{B}$ is an apartment containing both $z_1$ and $z_2$, then $I_{A'}(x) \subseteq I(x)$. Now $z_1, z_2 \in I_{A'}(x)$ and $I_{A'}(x)$ is convex. From this the third statement follows.

**3.15. Theorem.** *Let $x \in Y^{ss}$. Then $x \in Y^s$ if and only if $I(x)$ is bounded.*

*Proof.* We will first show that if $I(x)$ is not bounded, then $x \notin Y^s$. So let us assume that $I(x)$ is not bounded. Let $z \in I(x)$. We can find a sequence of points $z_i \in I(x)$ such that the distance dist $(z, z_i) \longrightarrow \infty$ if $i \longrightarrow \infty$. Let $A_i \subset \mathcal{B}$ be an apartment containing both $z$ and $z_i$.

Then $I_{A_i}(x) \subseteq I(x)$ and $z, z_i \in I_{A_i}(x)$. Let $A \subset \mathcal{B}$ be an apartment that contains the point $z$. There exist elements $h_i \in P_z$ such that $A_i = h_i A$. Since $P_z$ is compact, a subsequence of the sequence $h_i$ converges to an element $h \in P_z$. Clearly $I_{hA}(x) \subseteq I(x)$ and $I_{hA}(x)$ is not bounded. Hence $x \notin Y^s_{hA}$. So $x \notin Y^s$. This proves the "only if" part of the theorem.

Next we will show that if $x \in Y^{ss} - Y^s$, then the interval $I(x)$ is not bounded. Let $x \in Y^{ss} - Y^s$. There exists an apartment $\widetilde{A}$ such that $x \in Y^{ss}_{\widetilde{A}} - Y^s_{\widetilde{A}}$. In particular, there exists using lemma 3.17 below a 1-ps $\epsilon$ of the torus $\widetilde{S}$ that belongs to the apartment $\widetilde{A}$ such that $x \notin Y^s_{g\widetilde{A}}$ for all $g \in P(\epsilon)$. Let $z$ be a point in $I(x)$. Then there exists an element $g \in P(\epsilon)$ such that $z \in A' := g\widetilde{A}$. Since $x \in Y^{ss}$, we have $z \in Y^{ss}_{A'} - Y^s_{A'}$. Furthermore, $I_{A'}(x) \subseteq I(x)$, since $z \in A'$. Hence $I(x)$ is not bounded.

**3.16. Definition.** Let $S \subset G$ be a maximal $K^\circ$-split torus. For a 1-ps $\epsilon : (K^\circ)^* \to S$ we define the parabolic subgroup $P(\epsilon) := <U_\alpha(K), Z(K) \mid \alpha \in \Phi^+_K, <\alpha, \epsilon> \geq 0> \subset G(K)$.

Now we can state the following well-known lemma (See [MFK] appendix 2B, [K] theorem 4.2 and [He] theorem 5.2). For completeness we supply a proof.

**3.17. Lemma.** *Let $S \subset G$ be a maximal $K^\circ$-split torus and let $A \subset \mathcal{B}$ be the apartment that belongs to $S$. Let $x \in X \otimes_{L^\circ} \mathbb{K}$ be a point. Then the following two statements hold:*

i) *If $x \notin Y^s_A$, then there exists a 1-ps $\epsilon$ of $S$ such that $x \notin Y^s_{gA}$ for all $g \in P(\epsilon)$.*
ii) *If $x \notin Y^{ss}_A$, then there exists a 1-ps $\epsilon$ of $S$ such that $x \notin Y^{ss}_{gA}$ for all $g \in P(\epsilon)$.*

*Proof.* We only prove statement (i) for a point $x \in Y^{ss}_A - Y^s_A$. The polyhedron

$r(\mu_K(x))$ contains 0, but 0 is not contained in its interior. Hence one can find a 1-ps $\epsilon$ of $S$ such that $<\chi, \epsilon> \leq 0$ for all $\chi \in r(\mu_K(x)) \subset \mathcal{X}(S) \otimes \mathbb{R}$. We claim that this 1-ps satisfies statement (i) of the lemma.

The action of the additive groups $U_\alpha(K) \subset G(K)$ on the coordinates is given by: $u_\alpha(f_\alpha)^* x_{\beta,j} = x_{\beta,j} + \sum_\nu a_{\nu,i}(f_\alpha) x_{\nu,i}$, where $r(\nu) = r(\beta) - n\alpha$ for some $n \in \mathbb{Z}_{\geq 0}$ (See [Hu] proposition 27.2). Hence for an element $g \in P(\epsilon)$, one has that $<\chi, \epsilon> \leq 0$ for all $\chi \in r(\mu(g^{-1}(x)))$ and therefore $x \notin Y_{gA}^s$. Hence $x \notin Y_{gA}^s$ for all $g \in P(\epsilon)$.

### 3.18. Proposition.

i) A point $x \in Y^{ss}$ is contained in $Y_z^s$ if and only if $I(x) = \{z\}$.

ii) $Y_z^s \subset Y^s$.

iii) $Y^s = \bigcup_{z \in \mathcal{B}} Y_z^s$ if and only if $Y^s = Y^{ss}$.

*Proof.* A point $x \in Y^{ss}$ is in $Y_z^s$ if and only if for all apartments $A$ containing $z$, one has $I_A(x) = \{z\}$. This happens only when $I(x) = \{z\}$. This proves the first statement.

Let $x \in Y_z^s$. Then $I(x) = \{z\}$ and hence is bounded. By theorem 3.15 the point $x$ is in $Y^s$. Therefore $Y_z^s \subset Y^s$. This proves the second statement.

The third statement of the proposition follows from theorem 3.11 and the fact that if $Y^s \neq Y^{ss}$, then there exist points $x \in Y^s$ such that $I(x)$ does not consist of a single point.

**3.19. Remark.** Let $F \subset \mathcal{B}$ be a bounded subset. Let $cc(F)$ be the *circumcenter* of $F$. The circumcenter of $F$ is the center of the closed ball with minimal radius that contains $F$. It is well-defined and unique (See [B] Ch. VI §4). This gives a $G(K)$-equivariant map $\psi : Y^s \longrightarrow \mathcal{B}$, where $\psi(x) := cc(I(x))$. This map cannot be extended to all of $Y^{ss}$ if $Y^{ss} \neq Y^s$. This clarifies somewhat property 3.12 in [Ra].

For a point $z \in \mathcal{B}$ the analytical subspace $\psi^{-1}(z) \subset Y^s$ is not quasi-compact if $Y^s \neq Y^{ss}$. If $Y^s = Y^{ss}$, then $\psi^{-1}(z) = Y_z^s$ and therefore is quasi-compact. In [PV] and [V.1] it is shown that if $Y^s = Y^{ss}$, then the quotient $Y^s/\Gamma$ is a proper rigid analytic variety for all discrete co-compact subgroups $\Gamma \subset G(K)$.

**3.20. Remark.** One can define an interval of semistability for suitable subcomplexes $\mathcal{F} \subseteq \mathcal{B}$. Let $\mathcal{F} \subseteq \mathcal{B}$ be a subcomplex that satisfies the following conditions:

i) $\mathcal{F}$ is convex.

ii) $\mathcal{F} = \bigcup \{A \subset \mathcal{B} \mid A \subseteq \mathcal{F}\}$.

One puts $Y_{\mathcal{F}}^s := \bigcap_{A \subseteq \mathcal{F}} Y_A^s$ and $Y_{\mathcal{F}}^{ss} := \bigcap_{A \subseteq \mathcal{F}} Y_A^{ss}$. Furthermore, for a point $z \in \mathcal{F}$ one takes $Y_{z,\mathcal{F}}^s := \bigcap_{A \ni z,\, A \subseteq \mathcal{F}} Y_{z,A}^s$ and $Y_{z,\mathcal{F}}^{ss} := \bigcap_{A \ni z,\, A \subseteq \mathcal{F}} Y_{z,A}^{ss}$. Then $Y_{\mathcal{F}}^{ss} = \bigcup_{z \in \mathcal{F}} Y_{z,\mathcal{F}}^{ss}$ holds. Therefore one can define an interval of $\mathcal{F}$-semistability $I_{\mathcal{F}}(x)$ for $x \in Y_{\mathcal{F}}^{ss}$ as $I_{\mathcal{F}}(x) := \{z \in \mathcal{F} \mid x \in Y_{z,\mathcal{F}}^{ss}\}$. Then $I_{\mathcal{F}}(x)$ is non-empty and convex for $x \in Y_{\mathcal{F}}^{ss}$. Moreover, $I_{\mathcal{F}}(x)$ is bounded if and only if $x \in Y_{\mathcal{F}}^s$. The proofs given above for the case $\mathcal{F} = \mathcal{B}$ remain valid, mutatis mutandis.

**3.21. Compactification of the building.** In order to extend the map $I$ to the points of $X \otimes_{L^\circ} \mathbb{K}$ that are not contained in $Y^{ss}$, we compactify the building $\mathcal{B}$. This is done by adding the spherical building $\mathcal{B}_\infty$ to $\mathcal{B}$ (See [B] Ch. VI §9 or [Ro] Ch. 9 §3). Let us briefly recall the construction.

Let $S$ be a maximal $K^\circ$-split torus of $G$. Then $S$ determines an apartment $A \subset \mathcal{B}$ and a spherical apartment $A_\infty \subset \mathcal{B}_\infty$. We identify $A_\infty$ with $(\mathcal{X}_*(S) \otimes \mathbb{R} - \{0\})/\mathbb{R}_{>0}$. To a point $\delta \in A_\infty$, we associate the parabolic subgroup $P(\delta) := < Z(K), U_\alpha(K) \mid \alpha \in \Phi_K, < \alpha, \delta > \geq 0 >$ of $G(K)$. Here $Z$ denotes the centraliser of $S$ in $G$.

On the set of spherical apartments $g(A_\infty)$, $g \in G(K)$ one defines an equivalence relation as follows. Let $\delta_i \in g_i(A_\infty)$, $i = 1, 2$. Then $\delta_1 \sim \delta_2$ if and only if there exists an element $g \in P(\delta_1)$ such that $\delta_2 = g \cdot \delta_1 \cdot g^{-1}$. Note that $P(\delta_1) = P(\delta_2)$ if $\delta_1 \sim \delta_2$.

Now the *spherical building at infinity* is defined as $\mathcal{B}_\infty := \bigcup_{g \in G(K)} g(A_\infty)/\sim$. The compactification $\overline{\mathcal{B}}$ of $\mathcal{B}$ is obtained by adding $\mathcal{B}_\infty$ as a boundary to $\mathcal{B}$.

**3.22. Definition.** Two points $\delta_1, \delta_2 \in \mathcal{B}_\infty$ are called *antipodal* if there exists a spherical apartment containing both points and such that $\delta_1$ and $\delta_2$ are antipodal inside this apartment. A subset $F \subset \mathcal{B}_\infty$ is called *semi-convex* if it contains the line joining any pair of points in $F$ that are not antipodal. The *semi-convex hull* of a subset $F \subset \mathcal{B}_\infty$ is the smallest semi-convex subset of $\mathcal{B}_\infty$ that contains $F$.

For a point $x \in X \otimes_{L^\circ} \mathbb{K}$ we let $\widetilde{F}(x)$ denote the image of $\bigcup_{g \in G(K)} \widetilde{F}_{g(A)}(x)$ in $\mathcal{B}_\infty$. We define $\Lambda(x)$ as being the semi-convex hull of $\widetilde{F}(x)$. The map $\overline{I} : X \otimes_{L^\circ} \mathbb{K} \to$

{subsets of $\overline{\mathcal{B}}$} is defined as follows:

$$\overline{I}(x) := I(x) \cup \Lambda(x).$$

Here we take $I(x) := \emptyset$ for points $x \in X \otimes_{L^\circ} \mathbb{K} - Y^{ss}$.

**3.23. Proposition.** *Let $x \in X \otimes_{L^\circ} \mathbb{K}$ be a point. Then:*

i) $\Lambda(x) = \bigcup_{A \subset \mathcal{B}} \Lambda_A(x)$.

ii) $\Lambda(x) = \emptyset$ *if and only if* $x \in Y^s$.

iii) *If $x \in Y^{ss}$, then the points of $\Lambda(x)$ correspond to the rays that are contained in $I(x)$.*

*Proof.* This is a direct consequence of proposition 2.11.

**3.24. Remark.** The proposition above shows that $\overline{I}(x)$ is the closure in $\overline{\mathcal{B}}$ of $I(x)$ for points $x \in Y^{ss}$. However, the map $\overline{I}$ is in general not uniquely determined by the map $I$ (See remark 2.12).

## §4 Variation of line bundles

**4.1.** We study how the interval of semistability varies, when we change the line bundle $\mathcal{L}$. Our results in this case are direct translations of well-known results about variation of semistable points (See [DH] and [BP]).

**4.2. Equivalence classes of ample line bundles.** We call two ample line bundels $\mathcal{L}$ and $\mathcal{L}'$ *equivalent* if they determine the same variety of semistable points, that is, if $X^{ss}(S, \mathcal{L}) = X^{ss}(S, \mathcal{L}')$. If two line bundles are equivalent, then they also determine the same set of stable points.

Let $X = G/P_J$ with $J \subseteq \{1, \ldots, \ell\}$ a non-empty subset. Then the ample line $\mathcal{L}$ on $X$ corresponds to a unique weight $\lambda = \sum_{j \in J} n_j \cdot \omega_j$, with all $n_j > 0$ (See 1.4). The equivalence class of the line bundle $\mathcal{L}$ that corresponds to $\lambda$ will be denoted by $[\lambda]$. The ample line bundles that correspond to the weights $n \cdot \lambda$ with $n \in \mathbb{Z}_{>0}$ are all in the equivalence class $[\lambda]$.

By $C(X) \subset \mathcal{X}(T) \otimes \mathbb{Q}$ we denote the cone $\sum_{j \in J} n_j \cdot \omega_j$, where the $n_j$ are contained in $\mathbb{Q}_{>0}$. For an element $\lambda \in C(X)$ the equivalence class $[\lambda]$ will mean the equivalence class of the line bundle corresponding to the weight $n \cdot \lambda$, where $n \in \mathbb{Z}_{>0}$ is such that $n \cdot \lambda \in \mathcal{X}(T)$. The equivalence classes of ample line bundles on $X$ give a finite decomposition of $C(X)$. Indeed, it follows from propositions 1.7 and 1.9 that the equivalence class $[\lambda]$ is entirely determined by the set of hyperplanes spanned by roots that contain the weight $\lambda$.

In the remainder of this section we will write $X^s(S, [\lambda])$ and $X^{ss}(S, [\lambda])$ instead of $X^s(S, \mathcal{L})$ and $X^{ss}(S, \mathcal{L})$. Here $\lambda$ is the weight that corresponds to the line bundle $\mathcal{L}$. Furthermore we will write $\mu(x, \lambda)$ for the polyhedron $\mu(x) \subset \mathcal{X}(T) \otimes \mathbb{R}$ that is determined by the line bundle $\mathcal{L}$ corresponding to the weight $\lambda$ and a point $x \in X$.

The following proposition and its corollary are well-known:

**4.3. Proposition.** *Let $\lambda \in C(X) \subset \mathcal{X}(T) \otimes \mathbb{Q}$. Let $B(\lambda, \varepsilon) \subset \mathcal{X}(T) \otimes \mathbb{R}$ be the ball with centre $\lambda$ and radius $\varepsilon$. Then the following statements hold:*

i) *If $x \in X^s(S, [\lambda])$, then there exists an $\varepsilon > 0$ such that $x \in X^s(S, [\lambda'])$ for all $\lambda' \in B(\lambda, \varepsilon) \cap C(X)$.*

ii) *If $x \notin X^{ss}(S, [\lambda])$, then there exists an $\varepsilon > 0$ such that $x \notin X^{ss}(S, [\lambda'])$ for all $\lambda' \in B(\lambda, \varepsilon) \cap C(X)$.*

*Proof.* We may assume that $\lambda \in \mathcal{X}(T)$. The point $x$ is in $X^s(S, [\lambda])$ if and only if 0 is in the interior of the polyhedron $r(\mu(x, \lambda))$. The polyhedron $r(\mu(x, \lambda))$ is the convex hull of the weights $r(w(\lambda))$ such that $|x_{w(\lambda)}| = 1$. Now the polyhedron $r(\mu(x, \lambda'))$ is the convex hull of the weights $r(w(\lambda'))$, where the elements $w \in W$ occurring are the same as above. Hence if the difference between $\lambda$ and $\lambda'$ is very small, then the point 0 is still in the interior of $r(\mu(x, \lambda'))$. This proves (i).

The proof of part (ii) is similar. Now one uses the fact that a point $x$ is not in $X^{ss}(S, [\lambda])$ if and only if 0 is not contained in the polyhedron $r(\mu(x, \lambda))$. We leave the details to the reader.

**4.4. Corollary.** *Let $\lambda \in C(X) \subset \mathcal{X}(T) \otimes \mathbb{Q}$. Then there exists an $\varepsilon \in \mathbb{R}_{>0}$ such that for all $\lambda' \in B(\lambda, \varepsilon) \cap C(X)$ the following two statements hold:*

i) $X^s(S, [\lambda]) \subseteq X^s(S, [\lambda'])$

ii) $X^{ss}(S, [\lambda]) \supseteq X^{ss}(S, [\lambda'])$

*Proof.* This follows immediately from proposition 4.3 above and the fact that only a finite number of different polyhedra $r(\mu(x, \lambda))$ can occur (See also [DH] lemma

4.2.1).

**4.5. Definition.** On the set of equivalence classes of ample line bundles we put a partial ordering, which we denote by $\preceq$.

We write $[\lambda'] \preceq [\lambda]$ if $X^{ss}(S, [\lambda']) \supseteq X^{ss}(S, [\lambda])$. It follows from corollary 4.4 above that $[\lambda'] \preceq [\lambda]$ if and only if the area corresponding to $[\lambda']$ in $C(X)$ is in the closure of the area that corresponds to $[\lambda]$.

Till the end of this section we will write $Y^{ss}([\lambda])$ instead of $Y^{ss}$, if the line bundle $\mathcal{L}$ used to define $Y^{ss}$ is in the equivalence class $[\lambda]$. Similarly we write $I_A(x, [\lambda])$ and $I(x, [\lambda])$ instead of $I_A(x)$ and $I(x)$.

**4.6. Proposition.** *Let $x \in Y^{ss}([\lambda])$ and let $[\lambda'] \preceq [\lambda]$. Then one has:*

i) $I_A(x, [\lambda]) \subseteq I_A(x, [\lambda'])$

ii) $I(x, [\lambda]) \subseteq I(x, [\lambda'])$

*Proof.* Statement (i) follows directly from the definitions. Statement (ii) follows if one applies statement (i) to all apartments $A \subset \mathcal{B}$ such that the intersection $A \cap I(x, [\lambda'])$ is not empty.

**4.7. Theorem.** *Let $C(X, x) \subseteq C(X)$ denote the set of $\lambda \in C(X)$ such that $x \in Y^{ss}([\lambda])$. Then $\bigcup_{\lambda \in C(X,x)} I(x, [\lambda])$ is a connected subset of the building $\mathcal{B}$.*

*Proof.* This follows directly from the fact that $C(X, x)$ is convex and proposition 4.6(ii) above.

**4.8. Variation of homogeneous varieties.** Let $X := G/P_J$, $\emptyset \neq J \subseteq \{1, \ldots, \ell\}$. Let us assume that $\sharp J > 1$. Let $J' \subset J$, $J' \neq \emptyset, J$. We have a projection map $\pi_{J'} : X \longrightarrow X_{J'} := G/P_{J'}$.

Let $\mathcal{L}$ be an ample line bundle on $X_{J'}$. Then $\mathcal{L}$ corresponds to a weight $\lambda = \sum_{j \in J'} n_j \cdot \omega_j$, $n_j > 0$. The line bundle $\pi_{J'}^* \mathcal{L}$ on $X$ is not ample, but it is generated by global sections. It corresponds to the same weight $\lambda$ as $\mathcal{L}$ does.

Let $S \subset G$ be a maximal $K^\circ$-split torus. By $X^{ss}(S, [\lambda])$ and $X^s(S, [\lambda])$ we will denote the varieties $\pi_{J'}^{-1}(X_{J'}^{ss}(S, [\lambda]))$ and $\pi_{J'}^{-1}(X_{J'}^s(S, [\lambda]))$, respectively. These are the varieties of the points in $X$ that are semistable and stable, respectively, for the action of $S$ w.r.t. the line bundle $\pi_{J'}^* \mathcal{L}$.

One now defines the analytic subvarieties $Y^{ss}([\lambda])$ and $Y^{s}([\lambda])$ of $X \otimes_{L^\circ} \mathbb{K}$ that consist of the points that are semistable and stable, respectively, for all maximal $K$-split tori in $G \otimes_{K^\circ} K$. We will also denote the map $\pi_{J'} \otimes_{L^\circ} \mathbb{K} : X \otimes_{L^\circ} \mathbb{K} \longrightarrow X_{J'} \otimes_{L^\circ} \mathbb{K}$ by $\pi_{J'}$. One has $\pi_{J'}(Y^{ss}([\lambda])) = Y^{ss}(J', [\lambda])$ and $\pi_{J'}(Y^{s}([\lambda])) = Y^{s}(J', [\lambda])$. Here $Y^{ss}(J', [\lambda])$ and $Y^{s}(J', [\lambda])$ are the subspaces consisting of the points in $X_{J'} \otimes_{L^\circ} \mathbb{K}$ that are semistable and stable, respectively, for all maximal $K$-split tori w.r.t. the line bundle $\mathcal{L}$. For $x \in Y^{ss}([\lambda])$ one has also an interval of $G(K)$-semistability, defined by $I(x, [\lambda]) := I(\pi_{J'}(x), [\lambda])$.

Let $\overline{C}(X) \subset \mathcal{X}(T) \otimes \mathbb{Q}$ consist of the elements $\sum_{j \in J} n_j \cdot \omega_j$ with all $n_j \in \mathbb{Q}_{\geq 0}$ and with at least one $n_j \neq 0$. One verifies that all the results obtained above for ample line bundles (corresponding to weights $\lambda \in C(X)$) remain valid for line bundles on $X$ that correspond to weights $\lambda \in \overline{C}(X)$.

**4.9. Example.** Let $G$ be the split symplectic group $Sp_4$ and let $X = G/B$, where $B \subset G$ is a Borel subgroup. Let $P_1$ and $P_2$ be non-isomorphic maximal parabolic subgroups of $G$. Then $X_1 := G/P_1 \cong \mathbb{P}^3$ and $X_2 := G/P_2$ is given by a quadratic equation in $\mathbb{P}^4$. Let $\pi_i : X \longrightarrow X_i$, $i = 1, 2$ be the projection maps. All ample line bundles on $X$, $X_1$ and $X_2$ are in the single equivalence class $[\omega_1 + \omega_2]$, $[\omega_1]$ and $[\omega_2]$, respectively. Here $\omega_1$ and $\omega_2$ are the fundamental weights that correspond to the maximal parabolic subgroups $P_1$ and $P_2$, respectively.

One verifies that $X^{ss}(S, [\omega_1 + \omega_2]) = X^{s}(S, [\omega_1 + \omega_2])$. Furthermore, we have $X^{ss}(S, [\omega_1 + \omega_2]) = X^{ss}(S, [\omega_1]) \cap X^{ss}(S, [\omega_2])$. Therefore the interval of semistability $I(x, [\omega_1 + \omega_2])$ equals the intersection $I(x, [\omega_1]) \cap I(x, [\omega_2])$ for points $x \in Y^{ss}([\omega_1 + \omega_2])$. Moreover, this intersection consists of a single point, since $Y^{s}([\omega_1 + \omega_2]) = Y^{ss}([\omega_1 + \omega_2])$.

**4.10. Example.** Let $f(x, y)$ be a non-degenerate quadratic form in three variables that is defined over $K^\circ$. Let $L \supset K$ be a separable quadratic extension and let $\tau$ be the generator of $Gal(L/K)$. Then $f(x, \tau(y))$ is a non-degenerate unitary form.

Let $G$ be the unitary group that preserves the form $f(x, \tau(y))$ (considered as a group defined over $K^\circ$). Then $G \otimes_{K^\circ} L^\circ \cong SL_3$.

Let $X := G/B$, $X_1 := G/P_1 \cong \mathbb{P}^{\circ}_{L^\circ}$ and $X_2 := G/P_2 \cong \mathbb{P}^{2}_{L^\circ}$. Here $B \subset G$ is a Borel subgroup and $P_1$, $P_2 \subset G$ are maximal parabolic subgroups. Let $\omega_1$ and $\omega_2$ be the fundamental weights that correspond to the two maximal parabolic subgroups $P_1$ and $P_2$, respectively. The ample line bundles on $X$, $X_1$ and $X_2$ are all in the equivalence class $[\omega_1 + \omega_2]$, $[\omega_1]$ and $[\omega_2]$, respectively.

In this situation we have as in the previous example that:

$$Y^{ss}([\omega_1 + \omega_2]) = Y^s([\omega_1 + \omega_2]) = Y^{ss}([\omega_1]) \cap Y^{ss}([\omega_2]).$$

Therefore $I(x, [\omega_1 + \omega_2]) = I(x, [\omega_1]) \cap I(x, [\omega_2])$ for points $x \in Y^{ss}([\omega_1 + \omega_2])$. Again this intersection consists of a single point.

The variety $X = G/B \subset G/P_1 \times G/P_2 \cong \mathbb{P}^2 \times \mathbb{P}^2$ is given by the equation $f(x,y) = 0$. Here $x$ and $y$ are in the first and second $\mathbb{P}^2$, respectively. An element $g \in SU_3(L)$ acts on $\mathbb{P}^2_L \times \mathbb{P}^2_L$ as $(g, \tau(g))$. Let $\widetilde{\tau} \in Gal(K_s/K)$ be an element such that $\widetilde{\tau}|_L = \tau$.

Let $x$, $y \in X_1(K_s) \cong \mathbb{P}^2_L(K_s)$ be points that are semistable for all maximal $K$-split tori of $G \otimes_{K^\circ} K$. Then $(x, \widetilde{\tau}(y)) \in X(K_s)$ if and only if $f(x, \widetilde{\tau}(y)) = 0$. So it follows from $f(x, \widetilde{\tau}(y)) = 0$ that $I_1(x) \cap I_2(\widetilde{\tau}(y))$ consists of a single point. Here $I_1$ and $I_2$ are the intervals of semistability for the flag domains in $X_1$ and $X_2$, respectively. Since $g \in SU_3(L)$ acts on $X_1$ as $g$ and on $X_2$ as $\tau(g)$, one has $I_2(\widetilde{\tau}(y)) = I_1(y)$. It therefore follows from $f(x, \widetilde{\tau}(y)) = 0$, that $I_1(x) \cap I_1(y)$ consists of a single point (See also [LV] proposition 4.7).

## §5 Compactification

**5.1.** Let $P \subset G(K)$ be a parabolic subgroup. We call a maximal $K$-split torus $S_K \subset G \otimes_{K^\circ} K$ a *P-torus* if $S_K(K)$ is contained in $P$. An apartment $A \subset \mathcal{B}$ is called a *P-apartment* if the corresponding maximal $K$-split torus is a $P$-torus.

In this section we use the same notation as in §§1-3. In particular, we again fix an ample line bundle $\mathcal{L}$ on the projective variety $X$. We also fix again a field of definition for the pair $(X, \mathcal{L})$, which we denote by $L$.

The spaces $Y^s$ and $Y^{ss}$ consist of the points that are stable and semistable, respectively, for all maximal $K$-split tori in $G \otimes_{K^\circ} K$. However, it is sufficient to take the points in $X \otimes_{L^\circ} \mathbb{K}$ that are stable and semistable, respectively, for all maximal $K$-split $B$-tori for some fixed Borel subgroup $B \subset G(K)$ to obtain these spaces. This enables us to construct a $B$-equivariant compactification of $Y^s$.

Let $S \subset G$ be a maximal $K^\circ$-split torus and let $\chi \in \mathcal{X}(S)$ be a character. Let $P^\chi \subset G(K)$ be the parabolic subgroup that stabilises the character $\chi$. We translate the usual $S \otimes_{K^\circ} L^\circ$-linearization of $\mathcal{L}^{\otimes n}$, $n >> 0$ by the character $\chi$. Let $Y^s_\chi$ and $Y^{ss}_\chi$ be the subspaces of $X \otimes_{L^\circ} \mathbb{K}$ that consist of the points that are stable and semistable, respectively for all maximal $K$-split $P^\chi$-tori. Then $Y^s \subseteq Y^s_\chi \subseteq Y^{ss}_\chi \subseteq Y^{ss}$. There exist characters $\chi \in \mathcal{X}(S)$ such that $Y^s_\chi = Y^{ss}_\chi$. Then $Y^s_\chi$ is a $P^\chi$-equivariant compactification of $Y^s$.

We first show in propositions 5.2 and 5.3 below, that it suffices to use only $B$-apartments to define the analytic spaces $Y^{ss}$, $Y^s$, $Y^{ss}_z$ and $Y^s_z$. Here $B \subset G(K)$ is a Borel subgroup. Then we define and study the spaces $Y^{ss}_\chi$ and $Y^s_\chi$, where $\chi \in \mathcal{X}(S)$ is a character. Our methods here are very similar to those used in §2 and §3 to study the spaces $Y^{ss}$ and $Y^s$. In particular, we will also define and use an interval of semistability for these spaces.

**5.2. Proposition.** *Let $B \subset G(K)$ be a Borel subgroup and let $A \subset \mathcal{B}$ be a $B$-apartment. Then the following holds:*

i) $Y^s = \bigcap_{g \in B} g(Y^s_A)$.

ii) $Y^{ss} = \bigcap_{g \in B} g(Y^{ss}_A)$.

*Proof.* The proofs are in both cases very similar. Therefore we only prove case (i). Let $x \in X \otimes_{L^\circ} \mathbb{K} - Y^s$ be a point. Then there exists an apartment $A' \in \mathcal{B}$, such that $x \notin Y^s_{A'}$. Let $S'_K$ be the maximal $K$-split torus of $G \otimes_{K^\circ} K$ that corresponds to $A'$.

By lemma 3.17(i), there exists a 1-ps $\epsilon$ of $S'_K$ such that $x \notin g(Y^s_{A'})$ for all $g \in P(\epsilon)$. The parabolic subgroup $P(\epsilon)$ and $B$ have at least one maximal $K$-split torus in common. For this maximal $K$-split torus the point $x$ is not in the set of stable points. Hence $x \notin \bigcap_{g \in B} g(Y^s_A)$ and $Y^s \subseteq \bigcap_{g \in B} g(Y^s_A)$. Since the reverse inclusion is trivial, statement (i) of the proposition is true.

**5.3. Proposition.** *Let $B \subset G(K)$ be a Borel subgroup and let $A \subset \mathcal{B}$ be a $B$-apartment. Let $z \in A$ be a point. Then the following holds:*

i) $Y^s_z = \bigcap_{g \in P_z \cap B} g(Y^s_{z,A})$.

ii) $Y^{ss}_z = \bigcap_{g \in P_z \cap B} g(Y^{ss}_{z,A})$.

*Proof.* The proofs are in both cases very similar. Therefore we only prove case (i). Let $x \in X \otimes_{L^\circ} \mathbb{K} - Y^s_z$ be a point. Then there exists an apartment $A' \in \mathcal{B}$ that contains the point $z$, such that $x \notin Y^s_{z,A'}$. Let $S'_K$ be the maximal $K$-split torus of $G \otimes_{K^\circ} K$ that corresponds to $A'$.

We may assume that $0 \in A'$. Let $s \in S'_K(\mathbb{K})$ be such that $z = s \cdot 0$. Since $x \notin Y^s_{z,A'}$, there exists a 1-ps $\delta$ of $S'_K$ such that $<\chi, \delta> \leq 0$ for all $\chi \in r(\mu(s^{-1} \cdot x))$. Then $\delta(t) \cdot x \notin Y^s_{z,A'}$ for $\varepsilon < |t| \leq 1$, where $\varepsilon \in \mathbb{R}_{>0}$ is sufficiently close to 1.

Therefore $x \notin Y^s_{z',A'}$, where $z' = \delta(t)^{-1} \cdot z$ and $\varepsilon < |t| \leq 1$. We fix a point $z' \neq z$ as above, such that there exists a simplex $\sigma \in A'$ containing both $z$ and $z'$. There exists a $B$-apartment $A''$ that contains the simplex $\sigma$. Without loss of generality we may assume that $A = A''$. Let $g \in P_\sigma$ be such that $g(A') = A$.

Let $t_0 \in \mathbb{K}^*$ be such that $\delta(t_0)^{-1} \cdot z = z'$. Since $\delta(t)^{-1} \cdot z \in A' \cap g(A')$ for $|t_0| \leq |t| \leq 1$, the element $\delta(t) \cdot g \cdot \delta(t)^{-1}$ is contained in $s \cdot G(\mathbb{K}^\circ) \cdot s^{-1}$, where $s \in S'(\mathbb{K})$ is such that $s \cdot 0 = z$.

Therefore the reduction $\bar{g}$ of $g \in s \cdot G(\mathbb{K}^\circ) \cdot s^{-1}$ commutes with the 1-ps $\delta \otimes_{K^\circ} \overline{\mathbb{K}}$ of $S' \otimes_{K^\circ} \overline{\mathbb{K}}$. As a consequence $<\chi, \delta> \leq 0$ for all $\chi \in r(\mu(s^{-1} \cdot g^{-1}(x)))$. In particular, $x \notin Y^s_{z,g(A')} = Y^s_{z,A}$. Therefore $Y^s_z \subseteq \bigcap_{g \in P_z \cap B} g(Y^s_{z,A})$. Since the other inclusion is obvious, this proves part (i) of the proposition.

**5.4. Definition.** Let $S \subset G$ be a maximal $K^\circ$-split torus. Let $n \in \mathbb{Z}_{>0}$ and let $\chi \in \mathcal{X}(S)$. We consider the $S \otimes_{K^\circ} L^\circ$-linearization of the line bundle $\mathcal{L}^{\otimes n}$ obtained by translating the usual $S \otimes_{K^\circ} L^\circ$-linearization by the character $\chi$. Let $T \supseteq S$ be a maximal torus in $G$ that is defined over $K^\circ$. Then the linearization is such that $S \otimes_{K^\circ} L^\circ$ acts on the eigenspaces $V^\beta \subseteq V_{n \cdot \lambda}$, $\beta \in \mathcal{X}(T)$ for the usual action of $T \otimes_{K^\circ} L^\circ$ (coming form the $G \otimes_{K^\circ} L^\circ$ action) on $V_{n \cdot \lambda} = \oplus_{\beta \in \mathcal{X}(T)} V^\beta$ with the character $r(\beta) + \chi$.

The thus determined varieties of stable and semistable points in $X$ will be denoted by $X^s(S, \mathcal{L}, \chi/n)$ and $X^{ss}(S, \mathcal{L}, \chi/n)$, respectively. A point $x \in X$ is in $X^{ss}(S, \mathcal{L}, \chi/n)$ if and only if $-\chi/n$ is in $r(\mu(x))$. Moreover the point $x$ is stable if and only if $-\chi/n$ is in the interior of $r(\mu(x))$.

We will always assume that the integer $n > 0$ is large enough. To be more precise, we assume that $n > 0$ is such that the following two conditions hold:

i) $X^{ss}(S, \mathcal{L}, \chi/n) \subseteq X^{ss}(S, \mathcal{L})$.

ii) For all $m \in \mathbb{Z}$ such that $m \geq n$ one has $X^{ss}(S, \mathcal{L}, \chi/m) = X^{ss}(S, \mathcal{L}, \chi/n)$.

Applying corollary 4.4(ii) to this situation, one verifies that for any character $\chi \in$

$\mathcal{X}(S)$ there exists an integer $n > 0$ such that the two conditions above are satisfied.

With this restriction, our $S \otimes_{K^\circ} L^\circ$-linearization of the line bundle $\mathcal{L}^{\otimes n}$ is such that $X^{ss}(S, \mathcal{L}, \chi/n)$ only depends on the character $\chi$ and not on the value of $n \gg 0$. In particular, if the line bundle $\mathcal{L}$ is such that $X^s(S, \mathcal{L}) = X^{ss}(S, \mathcal{L})$, then $X^{ss}(S, \mathcal{L}, \chi/n) = X^{ss}(S, \mathcal{L})$ for $n \gg 0$.

Let $A \subset \mathcal{B}$ be the apartment that corresponds to the torus $S$. We write $Y^{ss}_{A,\chi}$ and $Y^s_{A,\chi}$ for the analytic spaces that correspond to $X^{ss}(S, \mathcal{L}, \chi/n) \otimes_{L^\circ} \mathbb{K}$ and $X^s(S, \mathcal{L}, \chi/n) \otimes_{L^\circ} \mathbb{K}$, respectively. Similarly, we write $Y^{ss}_{z,A,\chi}$ and $Y^s_{z,A,\chi}$ for the spaces that are the analogs in this situation of $Y^{ss}_{z,A}$ and $Y^s_{z,A}$.

**5.5. Lemma.** *Let $\chi \in \mathcal{X}(S)$ and let $z \in A$. Then the following holds:*

i) $X^s(S, \mathcal{L}) \subseteq X^s(S, \mathcal{L}, \chi/n) \subseteq X^{ss}(S, \mathcal{L}, \chi/n) \subseteq X^{ss}(S, \mathcal{L}, \chi/n)$ *and all inclusions are Zariski open.*

ii) $Y^s_{z,A} \subseteq Y^s_{z,A,\chi} \subseteq Y^{ss}_{z,A,\chi} \subseteq Y^{ss}_{z,A}$ *and all inclusions are open.*

iii) $Y^s_A \subseteq Y^s_{A,\chi} \subseteq Y^{ss}_{A,\chi} \subseteq Y^{ss}_A$.

iv) $Y^{ss}_{A,\chi} = \bigcup_{z \in A} Y^{ss}_{z,A,\chi}$.

*Proof.* The second inclusion in statement (i) is obvious. Moreover the third inclusion follows from our assumptions. The first inclusion again follows from the assumption that $X^{ss}(S, \mathcal{L}, \chi/n) \subseteq X^{ss}(S, \mathcal{L})$. Indeed, this is the analog in our situation of corollary 4.4(i). The openness of the inclusions is clear from the definitions. Statements (ii) and (iii) follow from (i).

The fourth statement of the proposition follows from the fact that the quotient $X^{ss}(S, \mathcal{L}, \chi/n)/(S \otimes_{K^\circ} L^\circ)$ is projective (See also prop. 2.4).

**5.6. Definition.** For $x \in Y^{ss}_A$ we put $n_A(x, \chi) := \sup \{<\chi, z> \mid z \in I_A(x)\}$. We also define an interval of semistability for $x \in Y^{ss}_{A,\chi}$ as being the subset $I_A(x, \chi) := \{z \in A \mid x \in Y^{ss}_{z,A,\chi}\}$ of $A$. It is clear that $I_A(x, \chi)$ is contained in $I_A(x)$ for any point $x \in Y^{ss}_{A,\chi}$. The following proposition and its corollary make this more precise.

**5.7. Proposition.** *Let $x \in Y^{ss}_{z,A}$ be a point. Then $x \in Y^{ss}_{z,A,\chi}$ if and only if $<\chi, z> = n_A(x, \chi)$.*

*Proof.* Without loss of generality we may assume that $z = 0$. Let us first consider the case where $x \in Y^{ss}_{0,A} - Y^{ss}_{0,A,\chi}$. Then $0 \in r(\mu(x))$ and $-\chi/m \notin r(\mu(x))$ for all $m \gg 0$. Therefore there exists a 1-ps $\epsilon$ of $S$, such that $<\delta, \epsilon> \geq 0$ for all characters $\delta \in r(\mu(x))$ and moreover $<\chi, \epsilon>\, > 0$. This implies that for $t \in \mathbb{K}^*$ with $\varepsilon < |t| \leq 1$, where $0 < \varepsilon < 1$, the point $\epsilon(t) \cdot x$ is in $Y_{0,A}^{ss}$. Hence $-\nu(\epsilon(t)) \in I_A(x)$. Since $<\chi, -\nu(\epsilon(t))>\, =\, v(\chi(\epsilon(t))) > 0$, we have $0 =<\chi, 0>\, <\, n_A(x, \chi)$.

Let us now consider the case, where the point $x$ is in $Y_{0,A,\chi}^{ss}$. Let $z \in I_A(x)$ be the point $s \cdot 0$ for some $s \in S(\mathbb{K})$. Then $s^{-1} \cdot x \in Y_{0,A}^{ss}$. Therefore $v(\delta(s^{-1})) \geq 0$ for all characters $\delta \in r(\mu(x))$. Since $-\chi/m \in r(\mu(x))$ for $m >> 0$, we have $v(\chi(s^{-1})) \leq 0$. Therefore $<\chi, z>\, =\, -v(\chi(s)) \leq 0$. In particular $0 =<\chi, 0>\, =\, n_A(x, \chi)$. This proves the proposition.

### 5.8. Corollary.

i) $Y_{A,\chi}^{ss} = \{x \in Y_A^{ss} \mid n_A(x, \chi) < \infty\}$.

ii) If $x \in Y_{A,\chi}^{ss}$, then $I_A(x, \chi) = \{z \in I_A(x) \mid <\chi, z>\, =\, n_A(x, \chi)\}$.

iii) Let $x \in Y_{A,\chi}^{ss}$. Then $I_A(x, \chi)$ is convex and non-empty.

*Proof.* This is a direct consequence of proposition 5.7.

### 5.9. Definition.
For a Borel subgroup $B \subset G(K)$ that contains $S(K)$ one defines $F_B$ as $F_B := \{\delta \in \mathcal{X}_*(S) \otimes \mathbb{R} \mid B \subseteq P(\delta)\}$. The image of $F_B$ in $A_\infty \subset \mathcal{B}_\infty$ is the simplex $\sigma_B \in \mathcal{B}_\infty$ that is stabilised by the Borel group $B$. To a character $\chi \in \mathcal{X}(S)$ we associate a parabolic subgroup $P^\chi \subset G(K)$. Let $P^\chi :=<B \mid \forall (\delta \in F_B)\;<\chi, \delta>\, \geq 0>$. In fact there exists a $\delta \in A_\infty \cong (\mathcal{X}_*(S) \otimes \mathbb{R} - \{0\})/\mathbb{R}_{>0}$ such that $P^\chi = P(\delta)$. Indeed, the intersection of the simplices $\sigma_B \in \mathcal{B}_\infty$ with $B \subseteq P^\chi$ is non-empty. Let us denote this intersection by $\tau_\chi$. Then $\tau_\chi \in \mathcal{B}_\infty$ is the simplex stabilised by $P^\chi$. Hence for a sufficiently general point $\delta \in \tau_\chi$ one has that $P(\delta) = P^\chi$.

For a simplex $\sigma \in \mathcal{B}_\infty$ and a point $z \in \mathcal{B}$ we denote by $C(z, \sigma)$ the cone in $\mathcal{B}$ that consists of the halflines that start in $z$ and end in $\sigma$. Here we say that a halfline ends in $\sigma$, if the corresponding equivalence class of parallel halflines is a point of $\sigma$.

### 5.10. Lemma.
Let $z \in \mathcal{B}$ be a point that is not contained in $A$ and let $z_A \in A$ be the point closest to $z$. Let $H \subset P^\chi \cap P_{z_A}$ be the subset that consists of the elements $g$ such that $g(z) \in A$. Then the value of $<\chi, g(z)>$ does not depend on the element $g \in H$.

*Proof.* Since $g \in P^\chi$ fixes $z_A$, the element $g$ fixes the cone $C(z_A, \tau_\chi)$. Furthermore, the cone $g(C(z, \tau_\chi)) = C(g(z), \tau_\chi)$ is contained in $A$.

If $g_1, g_2 \in H$, then $C(g_i(z), \tau_\chi) \subset A$ and $g_1 g_2^{-1} \in P^\chi \cap P_{z_A}$ maps $C(g_2(z), \tau_\chi)$ to $C(g_1(z), \tau_\chi)$. There exists an element $h \in P^\chi \cap P_{z_A}$, such that $hg_1 g_2^{-1}$ preserves $A$ and maps $C(g_2(z), \tau_\chi)$ to $C(g_1(z), \tau_\chi)$. Hence $hg_1 g_2^{-1} \in N(S) \cap P^\chi \cap P_{z_A}$. Since the intersection of the Weyl group $W_K$ of $S$ with $P^\chi$ preserves the character $\chi$, we have $<\chi, g_1(z)> = <\chi, g_2(z)>$. This proves the lemma.

**5.11. Definition.** Let $z \in \mathcal{B}$ be a point that is not contained in $A$ and let $z_A \in A$ be the point in $A$ closest to $z$. Then we put $f_\chi(z) := <\chi, g(z)>$, where $g$ is an element of $H$. Here $H$ is as in the lemma above. If $z \in A$, then we put $f_\chi(z) := <\chi, z>$. It follows from the lemma above that the function $f_\chi$ is well-defined.

We fix a $G(K)$-invariant distance $d(-,-)$ on the building $\mathcal{B}$. For two apartments $A, A' \subset \mathcal{B}$ we put $d(A, A') := \inf \{d(z, z') \mid z \in A, \ z' \in A'\}$.

In the following three propositions we study the function $f_\chi$ and its relation with $P^\chi$.

**5.12. Proposition.** *Let $g \in P^\chi$ and let $g(A)$ be a $P^\chi$-apartment. Then $f_\chi(z) = <\chi, g^{-1}(z)> + c_g$. Here $c_g$ is a constant depending on $g$.*

*Proof.* Let us first consider the case where the intersection $A \cap g(A)$ is non-empty. Then there exists an element $h \in G(K)$ that stabilises $A \cap g(A)$ and such that $h(g(A)) = A$. Since both $A$ and $g(A)$ are $P^\chi$-apartments, the element $h \in G(K)$ is actually contained in $P^\chi$.

Clearly $f_\chi(z) = <\chi, h(z)>$ for all $z \in g(A)$. Furthermore, there exists an element $h' \in N(S) \cap P^\chi$, such that $g(A \cap g(A)) = h'(h(A \cap g(A)))$. Then $h' = w \cdot s$ with $s \in S(\mathbb{K})$ and $w \in N(S) \cap P^\chi$ stabilises some point in $h(A \cap g(A))$. Therefore $f_\chi(z) = <\chi, g^{-1}(z)> - <\chi, \nu(s)>$ and the lemma holds in this case.

Let us now assume that $A \cap g(A) = \emptyset$. Let $F_A \subset A$ consist of the points $z \in A$ such that $d(z, g(A)) = d(A, g(A))$. Similarly, we denote by $F_{g(A)}$ the set of points $z \in g(A)$ such that $d(A, z) = d(A, g(A))$. The sets $F_A$ and $F_{g(A)}$ are both convex. Moreover, since $A \cap g(A) = \emptyset$, both $F_A$ anf $F_{g(A)}$ do not contain a sector. In particular, the parabolic subgroup $P^\chi$ is not a Borel subgroup.

Let $H_A \subset A$ be a wall that contains $F_A$. We can choose the wall $H_A$ in such a way that the reflection with respect to $H_A$ corresponds to an element of $N(S) \cap P^\chi$. Since $F_{g(A)}$ is parallel to $F_A$, there exists a wall $H_{g(A)} \subset g(A)$ of the same type as $H_A$ that contains $F_{g(A)}$.

Let $H^+_{g(A)}$ and $H^-_{g(A)}$ be the two halfspaces in $g(A)$ determined by the wall $H_{g(A)}$. Then there exists elements $h^+$ and $h^-$ in $G(K)$ that preserve $F_A$ and such that $h^+(H^+_{g(A)}) \subset A$ and $h^-(H^-_{g(A)}) \subset A$. Since the elements $h^+$ and $h^-$ stabilise $F_A$, they are contained in $P^\chi$. Then $f_\chi(z) = <\chi, h^+(z)>$ for $z \in H^+_{g(A)}$ and $f_\chi(z) = <\chi, h^-(z)>$ for $z \in H^-_{g(A)}$. One can moreover choose the elements $h^+$ and $h^-$ in such a way that $h^+(H_{g(A)}) = h^-(H_{g(A)})$. Since $g(H_{g(A)}) = w(h^-(H_{g(A)})) = w(h^+(H_{g(A)}))$ for some element $w \in N(S) \cap P^\chi$, the proposition follows.

**5.13. Remark.** On the set of cones $C(z, \sigma)$, $z \in \mathcal{B}$ one can define an equivalence relation. Take $C(z_1, \sigma) \sim C(z_2, \sigma)$ if and only if the intersection $C(z_1, \sigma) \cap C(z_2, \sigma)$ is non-empty. If a point $z_3$ is in the intersection, then the cone $C(z_3, \sigma)$ is contained in the intersection. The parabolic group $P(\sigma)$ acts on the set $\{C(z, \sigma) \mid z \in \mathcal{B}\}$. In fact, $\{C(z, \sigma) \mid z \in \mathcal{B}\}/\sim$ is the affine building of $P(\sigma)/R(P(\sigma))$. Here $R(P(\sigma))$ is the radical of the parabolic subgroup $P(\sigma) \subset G(K)$. This construction is described in [Ro] Ch. 10 §2 for the case where $\sigma$ is a codimension one simplex of $\mathcal{B}_\infty$. The general case is treated in [La] §2 and §13, where it is used to compactify the building $\mathcal{B}$.

The following proposition clarifies the relation between $f_\chi$ and $P^\chi$ somewhat. We omit the proof, since we will not need the result in the sequel. Furthermore, the proof is in the same vein as that of the previous proposition.

**5.14. Proposition.** *Let $g \in P^\chi$. Then*

i) *There exists a constant $c_g$ such that $f_\chi(g(z)) = f_\chi(z) - c_g$ for all points $z \in \mathcal{B}$.*

ii) *If $g$ is contained in a Levi subgroup of $P^\chi$, then $f_\chi(g(z)) = f_\chi(z)$ for all $z \in \mathcal{B}$.*

**5.15. Proposition.** *The subset $\{z \in \mathcal{B} \mid f_\chi(z) \geq n\}$ is convex for all $n \in \mathbb{R}$.*

*Proof.* Let $z_1, z_2 \in \mathcal{B}$ be two points such that $f_\chi(z_i) \geq n$ for $i = 1, 2$. To prove the proposition, it is sufficient to show that for every point $z$ in the shortest path $[z_1, z_2]$ in the building that joins the points $z_1$ and $z_2$ one has $f_\chi(z) \geq n$. We assume that $f_\chi(z_1) \geq f_\chi(z_2) \geq n$.

Without loss of generality, we may assume that $z_1$ is contained in the $P^\chi$-

apartment A. Let us first assume that $C(z_1, \tau_\chi) \cap [z_1, z_2] = \{z_1\}$. Then there exists a sector $C \subset A$ that starts in $z_1$, contains the cone $C(z_1, \tau_\chi)$ and such that $C \cap [z_1, z_2] = \{z_1\}$. Then there exists an element $h \in G(K)$ such that $h(C) = C$ and such that $h([z_1, z_2]) \subset A$. Since $h$ stabilises $C$, the element $h$ is contained in $P^\chi$. Therefore $h^{-1}(A)$ is a $P^\chi$-apartment that contains the path $[z_1, z_2]$. Using proposition 5.12, it is clear that for every point $z \in [z_1, z_2]$ one has $f_\chi(z_1) \geq f_\chi(z) \geq f_\chi(z_2) \geq n$.

Let us now assume that $C(z_1, \tau_\chi) \cap [z_1, z_2] = [z_1, z_3]$ for some point $z_3 \neq z_1$. Then $[z_1, z_3]$ is contained in the apartment $A$, since $A$ is a $P^\chi$-apartment. Moreover, $f_\chi(z_3) > f_\chi(z_1)$, since $z_3 \in C(z_1, \tau_\chi)$. In particular, $f_\chi(z_3) \geq f_\chi(z) \geq f_\chi(z_1) \geq n$ for all points $z \in [z_1, z_3]$.

Furthermore, $C(z_3, \tau_\chi) \cap [z_1, z_2] = \{z_3\}$. Therefore there exists a $P^\chi$-apartment $A'$ that contains the path $[z_3, z_2]$. Then for $z \in [z_3, z_2]$ one has $f_\chi(z_3) \geq f_\chi(z) \geq f_\chi(z_2) \geq n$. So for all points $z$ in the path $[z_1, z_2]$, we have $f_\chi(z) \geq n$. This proves the proposition.

**5.16. Definition.** Let $A \subset \mathcal{B}$ be a $P^\chi$-apartment. If $g_1, g_2 \in P^\chi$ are such that $g_1(A) = g_2(A)$, then $g_1(Y^{ss}_{A,\chi}) = g_2(Y^{ss}_{A,\chi})$. We define $Y^{ss}_\chi := \bigcap_{g \in P^\chi} g(Y^{ss}_{A,\chi})$ and $Y^s_\chi := \bigcap_{g \in P^\chi} g(Y^s_{A,\chi})$. Furthermore, we define $n(x, \chi) := \sup \{f_\chi(z) \mid z \in I(x)\}$ for points $x \in Y^{ss}$.

**5.17. Proposition.** $Y^{ss}_\chi = \{x \in Y^{ss} \mid n(x, \chi) < \infty\}$.

*Proof.* Let $x \in Y^{ss}$ be a point such that $x \notin Y^{ss}_\chi$. We will show that $n(x, \chi) = \infty$. Since $x \notin Y^{ss}_\chi$, there exists a $P^\chi$-apartment $A \subset \mathcal{B}$, such that $x \notin Y^{ss}_{A,\chi}$. In particular, $\sup \{f_\chi(z) \mid z \in I_A(x)\} = \infty$. There exist a point $z_0 \in I_A(x)$ and a 1-ps $\epsilon$ of the torus $S$ that belongs to $A$ such that the halfline $L_\epsilon := \{\epsilon(t) \cdot z_0 \mid t \in \mathbb{K}^*, |t| \leq 1\}$ is contained in $I_A(x)$. The 1-ps $\epsilon$ can be choosen in such a way that $\sup \{f_\chi(z) \mid z \in L_\epsilon\} = \infty$.

Since $L_\epsilon \subseteq I_A(x)$, we have $< \chi', \epsilon > \leq 0$ for all characters $\chi' \in r(\mu_K(x))$. By proposition 3.17(i), $x \in Y^{ss}_{A'} - Y^s_{A'}$ for all $P(\epsilon)$-apartments $A' \subset \mathcal{B}$.

Let $z_1 \in I(x)$ be a point. There exists a $P(\epsilon)$-apartment $A''$ that contains $z_1$. Let $g \in P(\epsilon)$ be such that $A'' = g(A)$ (Note that $A$ is also a $P(\epsilon)$-apartment). Let $\epsilon'$ be the 1-ps $g \circ \epsilon \circ g^{-1}$ of the torus $g \cdot S \cdot g^{-1}$. Let $L_{\epsilon'}$ be the halfline $L_{\epsilon'} := \{\epsilon'(t) \cdot z_1 \mid t \in$ $\mathbb{K}^*$, $|t| \leq 1\}$. Then $L_{\epsilon'} \subseteq I(x)$. Moreover, $sup\ \{f_\chi(z) \mid z \in L_{\epsilon'}\} = \infty$. Therefore $n(x, \chi) = \infty$.

To prove the reverse inclusion, we consider a point $x \in Y^{ss}$ with $n(x, \chi) = \infty$. Let $z_i \in I(x)$, $i \in \mathbb{N}$ be a sequence of points, such that $f_\chi(z_i) \longrightarrow \infty$ for $i \longrightarrow \infty$. Let $A_i \subset \mathcal{B}$ be an apartment that contains the points $z_1$ and $z_i$. For $i > 1$ there exists an element $h_i \in P_{z_1}$ such that $A_i = h_i(A_1)$. Since $P_{z_1}$ is compact, a subsequence of the sequence $h_i$, $i > 1$ converges to an element $h \in P_{z_1}$.

The apartment $h(A_1)$ contains infinitely many of the points $z_i$. Therefore one has: $sup\ \{f_\chi(z) \mid z \in I_{h(A_1)}(x)\} = \infty$. There exists a halfline $H \subseteq I_{h(A_1)}(x) \subseteq I(x)$, such that $sup\ \{f_\chi(z) \mid z \in H\} = \infty$. Some halfline $H' \subseteq H$ is contained in a $P^\chi$-apartment $\widetilde{A}$. Therefore $x \notin Y^{ss}_{\widetilde{A},\chi}$. So $x \notin Y^{ss}_\chi$.

**5.18. Definition.** Let $z \in \mathcal{B}$ be a point and let $A \subset \mathcal{B}$ be a $P^\chi$-apartment that contains $z$. Then we define:

$$Y^{ss}_{z,\chi} := \bigcap_{g \in P_z \cap P^\chi} g(Y^{ss}_{z,A,\chi}).$$

$$Y^s_{z,\chi} := \bigcap_{g \in P_z \cap P^\chi} g(Y^s_{z,A,\chi}).$$

We also define an interval of $P^\chi$-semistability $I(x, \chi)$ as follows:

$$I(x, \chi) := \{z \in \mathcal{B} \mid x \in Y^{ss}_{z,\chi}\}.$$

**5.19. Proposition.**

i) $Y^s_z \subseteq Y^s_{z,\chi} \subseteq Y^{ss}_{z,\chi} \subseteq Y^{ss}_z$ and all inclusions are open.

ii) Let $x \in Y^{ss}_z$ be a point. Then $x \in Y^{ss}_{z,\chi}$ if and only if $f_\chi(z) = n(x, \chi)$.

*Proof.* Statement (i) is a direct consequence of proposition 5.3 and lemma 5.5. So let us consider the second statement.

Let $x \in Y^{ss}_z$. If $f_\chi(z) = n(x, \chi)$, then $x \in Y^{ss}_{z,A,\chi}$ for all $P^\chi$-apartments $A$ that contain $z$. Hence $x \in Y^{ss}_{z,\chi}$. This proves one direction of the equivalence.

Let us now assume that $f_\chi(z) < n(x, \chi)$. Then there exists a point $z' \in I(x)$ such that $f_\chi(z') > f_\chi(z)$. Let $[z, z']$ be the shortest path joining the points $z$ and $z'$. Since $I(x)$ is convex, the path $[z, z']$ is contained in $I(x)$.

If the path $[z, z']$ is contained in some $P^\chi$-apartment $A$, then $f_\chi(z) < n_A(x, \chi)$. In particular, $x \notin Y^{ss}_{z,A,\chi}$.

So let us now assume that the path $[z, z']$ is not contained in any $P^\chi$-apartment $A'$. It follows from the proof of proposition 5.15, that $C(z, \tau_\chi) \cap [z, z'] \neq \{z\}$. In particular, there exists a point $z'' \neq z$ that is contained in $C(z, \tau_\chi) \cap [z, z']$. Let $A \subset \mathcal{B}$ be a $P^\chi$-apartment that contains the point $z$. Then $z'' \in A$. Furthermore, $f_\chi(z'') > f_\chi(z)$. Therefore $f_\chi(z) < n_A(x, \chi)$ and $x \notin Y^{ss}_{z,A,\chi}$.

We have now shown that if $x \in Y^{ss}_z$ and $f_\chi(z) < n(x, \chi)$, then $x \notin Y^{ss}_{z,A,\chi}$ for some $P^\chi$-apartment $A$ that contains $z$. Therefore $x \notin Y^{ss}_{z,\chi}$.

**5.20. Theorem.** $Y^{ss}_\chi = \bigcup_{z \in \mathcal{B}} Y^{ss}_{z,\chi}$.

*Proof.* Let $x \in Y^{ss}_\chi$ be a point. Then $n(x, \chi) < \infty$ by proposition 5.17. There exists a point $z \in I(x)$ such that $f_\chi(z) = n(x, \chi)$. Hence $x \in Y^{ss}_{z,\chi}$. So $Y^{ss}_\chi \subseteq \bigcup_{z \in \mathcal{B}} Y^{ss}_{z,\chi}$.

Let us now assume that $x \in Y^{ss}_{z,\chi}$ for some point $z \in \mathcal{B}$. Then $x \in Y^{ss}_z \subset Y^{ss}$ and $z \in I(x)$ is such that $f_\chi(z) = n(x, \chi)$. In particular, $n(x, \chi) < \infty$ and $x \in Y^{ss}_\chi$. Hence $\bigcup_{z \in \mathcal{B}} Y^{ss}_{z,\chi} \subseteq Y^{ss}_\chi$. This proves the theorem.

**5.21. Proposition.** *Let $x \in Y^{ss}_\chi$ be a point. Then:*

i) $I(x, \chi) = \{z \in I(x) \mid f_\chi(z) = n(x, \chi)\}$.

ii) $I(x, \chi)$ *is convex and non-empty.*

*Proof.* The first statement of the proposition is a direct consequence of propostion 5.19(ii). So let us consider statement (ii).

Using theorem 5.20, one shows that $I(x, \chi)$ is non-empty for $x \in Y^{ss}_\chi$. By proposition 5.19, $I(x, \chi)$ is the intersection of $I(x)$ and the set $F_\chi(x) := \{z \in \mathcal{B} \mid f_\chi(z) \geq n(x, \chi)\}$. Since $I(x)$ is convex by theorem 3.14(iii) and $F_\chi(x)$ is convex by proposition 5.15, the intersection $I(x, \chi)$ is convex.

**5.22. Proposition.** *There exist characters $\chi \in \mathcal{X}(S)$ such that $Y^s_\chi = Y^{ss}_\chi$.*

*Proof.* Let $x \in Y^{ss}$ be a point. Let $\chi \in \mathcal{X}(S)$ be a character. By proposition 2.7, the interval of $S$-semistability $I_A(x)$ is bounded by hyperplanes parallel to walls. Therefore, if the character $\chi$ is not contained in a hyperplane in $\mathcal{X}(S) \otimes \mathbb{R}$ that is spanned by roots, then $I_A(x, \chi)$ consists of a single point for $x \in Y^{ss}_{A,\chi}$. In particular, $Y^s_{A,\chi} = Y^{ss}_{A,\chi}$. Then also $Y^s_\chi = Y^{ss}_\chi$.

**5.23. Theorem.** *Let $\chi \in \mathcal{X}(S)$ be a character such that $Y_\chi^s = Y_\chi^{ss}$. Then there exists a formal scheme $\mathcal{Y}_\chi$ over $spf(\mathbb{K}^\circ)$ that has the following properties:*

i) *The generic fibre $\mathcal{Y}_\chi \otimes_{\mathbb{K}^\circ} \mathbb{K}$ of $\mathcal{Y}_\chi$ is $Y_\chi^s$.*

ii) *The closed fibre $\mathcal{Y}_\chi \otimes_{\mathbb{K}^\circ} \overline{\mathbb{K}}$ of $\mathcal{Y}_\chi$ consists of proper components that are in 1-1 correspondance with the vertices of the building $\mathcal{B}$ of $G(K)$.*

iii) *The parabolic subgroup $P^\chi \subset G(K)$ acts on $\mathcal{Y}_\chi$.*

*Proof.* We prove the theorem by constructing an admissable and pure affinoid covering of $Y_\chi^s$ that is $P^\chi$-invariant and that is such that the reduction has property (ii) of the theorem. We will follow [PV] sections 3.3-3.6 closely.

First we describe a pure affinoid covering of $Y_{A,\chi}^s$ for a $P^\chi$-apartment $A \subset \mathcal{B}$. For a simplex $\sigma \in A$ we take $Y_{\sigma,A,\chi} := \bigcup_{z \in \sigma} Y_{z,A,\chi}^s$. Then $Y_{A,\chi}^s = \bigcup_{\sigma \in A} Y_{\sigma,A,\chi}$. Furthermore, $Y_{\sigma_1,A,\chi} \cap Y_{\sigma_2,A,\chi}$ is non-empty if and only if $\sigma_1 \cap \sigma_2 \neq \emptyset$. If the intersection is non-empty, then the intersection equals $Y_{\sigma_1 \cap \sigma_2, A, \chi}$.

As in [PV] section 3.4, one shows that one can refine this covering of $Y_{A,\chi}^s$ into a pure affinoid covering of $Y_{A,\chi}^s$ in such a way that the reduction consists of proper components that are in a 1-1 correspondence with the vertices of the apartment $A$. The components of the reduction of $Y_{\sigma,A,\chi}$ correspond to the vertices that are contained in the simplex $\sigma$.

The next step in the proof is to use the spaces $Y_{\sigma,A,\chi}$ to construct a pure affinoid covering of $Y_\chi^s$ that has the desired properties. To a simplex $\sigma \in \mathcal{B}$ we associate the subspace $Y_{\sigma,\chi} := \bigcup_{z \in \sigma} Y_{z,\chi}^s$ of $Y_\chi^s$. Let $A \subset \mathcal{B}$ be a $P^\chi$-apartment that contains $\sigma$. Then $Y_{\sigma,\chi} = \bigcap_{g \in P_\sigma \cap P^\chi} g(Y_{\sigma,A,\chi}) - \bigcup_{\tau \subset \sigma} (\bigcap_{g \in P_\sigma \cap P^\chi} g(Y_{\tau,A,\chi}) - \bigcap_{g \in P_\tau \cap P^\chi} g(Y_{\tau,A,\chi}))$. Here the $\tau \subset \sigma$ are simplices that are contained in $\sigma$.

Then $Y_{\sigma,\chi} \subset Y_{\sigma,A,\chi}$ is an open subset. The covering $Y_{\sigma,\chi}$, for simplices $\sigma \in \mathcal{B}$ can be refined into a pure affinoid covering of $Y_\chi^s$ such that the components of the reduction correspond 1-1 to the vertices of the building. Moreover, the affinoid covering can be taken to be $P^\chi$-invariant.

To prove that the components of the reduction are proper, one proceeds as in [PV] section 3.6(g). One uses the fact that $Y_\chi^s = \bigcup_{\sigma \in \mathcal{B}} Y_{\sigma,\chi}$ equals $X \otimes_{L^\circ} \mathbb{K}$ minus a compact family of Zariski-closed subsets. From this it follows, that there exist admissible affinoid coverings $\{F_i\}$ and $\{H_i\}$ of $Y_\chi^s$, such that $F_i \subset\subset H_i$ for all $i$.

Then the properness of the components of the reduction follows from [Lü].

**5.24. Remark.** If $\chi \in \mathcal{X}(S)$ is a character such that $Y_\chi^s \neq Y_\chi^{ss}$, then the construction in the proof of the theorem above still gives a formal scheme $\mathcal{Y}_\chi$. The generic fibre of the formal scheme $\mathcal{Y}_\chi$ is $\bigcup_{z \in \mathcal{B}} Y_{z,\chi}^s$. So the generic fibre consists of the points $x \in Y_\chi^s$ such that $I(x, \chi)$ consists of a single point. The closed fibre still consists of components that are in 1-1 correspondance with the vertices of the affine building $\mathcal{B}$. These components are not proper anymore.

If $Y^s \neq Y^{ss}$, then the construction in [PV] section 3.6 still gives a formal scheme $\mathcal{Y}$. Its generic fibre is $\bigcup_{z \in \mathcal{B}} Y_z^s$. Therefore the generic fibre consists of the points $x \in Y^s$ for which $I(x)$ consists of a single point. The closed fibre again consists of components that are in 1-1 correspondance to the vertices of the building $\mathcal{B}$. These components are not proper, if $Y^s \neq Y^{ss}$.

The formal scheme $\mathcal{Y}$ is an open subscheme of the formal scheme $\mathcal{Y}_\chi$ for all characters $\chi \in \mathcal{X}(S)$. This justifies calling $Y_\chi^s$ a compactification of $Y$, if $\chi \in \mathcal{X}(S)$ is such that $Y_\chi^s = Y_\chi^{ss}$.

It is tempting to compare theorem 5.23 with the real case. In particular, in the case of the Satake and Baily-Borel compactifications of quotients of hermitian symmetric spaces by arithmetic subgroups, one adds for each parabolic subgroup of $G(\mathbb{Q})$ a boundary component to the hermitian symmetric space. In our situation it is not clear wether one can somehow glue the spaces $g(Y_\chi^s)$, $g \in G(K)$ together in some satisfactory way. However, for any character $\chi \in \mathcal{X}(S)$ the following two statements hold:

i) $Y^s = \bigcap_{g \in G(K)} g(Y_\chi^s) = \bigcap_{g \in G(K)} g(Y_\chi^{ss})$.

ii) If $S \otimes_{K^\circ} L^\circ$ acts without fixed points on $X^{ss}(S, \mathcal{L})$, then $Y^{ss} = \bigcup_{g \in G(K)} g(Y_\chi^{ss})$.

It is quite likely, that $Y^{ss} = \bigcup_{g \in G(K)} g(Y_\chi^{ss})$ remains true if $X^{ss}(S, \mathcal{L})$ does contain fixed points for the action of $S \otimes_{K^\circ} L^\circ$.

**REFERENCES**

[BGR] S. Bosch, U. Güntzer and R. Remmert, Non-archimedean analysis, Springer

Verlag, 1984.

[BP] M. Brion and C. Procesi, Action d'un tore dans une variété projective, in "Operator algebras, unitary representations, enveloping algebras and invariant theory", Birkhäuser 192 (1990), 509-539.

[B] K.S. Brown, Buildings, Springer Verlag, 1989.

[BrT] F. Bruhat and J. Tits, Groupes réductifs sur un corps local I : Données radicielles valuées, Publ. Math. I.H.E.S. 41 (1972), 5-121.

[DH] I.V. Dolgachev and Y. Hu, Variation of geometric invariant theory quotients, Publ. Math. I.H.E.S. 87 (1998), 5-56.

[D] V.G. Drinfel'd, Elliptic modules, Math. USSR Sbornik 23 (1974), 561-592.

[FH] H. Flaschka and L. Haine, Torus orbits in G/P, Pac. J. of math. 149 (1991), 251-292.

[FP] J. Fresnel and M. van der Put, Geométrie analytique rigide et applications, Progres in Math. 10, Birkhäuser, 1981.

[GS] I.M. Gel'fand and V.V. Serganova, Combinatorical geometries and torus strata on homogeneous compact manifolds, Russian Math. Surveys 42 (1987), 133-168.

[He] W.H. Hesselink, Uniform instability in reductive groups, J. f. d. reine u. angew. Math. 303/304 (1978), 74-96.

[Hu] J.E. Humphreys, Linear algebraic groups, Springer Verlag, 1975.

[J] J.C. Jantzen, Representations of algebraic groups, Academic Press, 1987.

[K] G.R. Kempf, Instability in invariant theory, Ann. of Math. 108 (1978), 299-316.

[LV] K.F. Lai and H. Voskuil, P-adic automorphic functions for the unitary group in three variables, Algebra Colloquium 7 (2000), 335-360.

[La] E. Landvogt, A compactification of the Bruhat-Tits building, Lect. Notes in Math. 1619 (1996).

[Lü] W. Lütkebohmert, Formal algebraic and rigid analytic geometry, Math. Ann. 286 (1990), 341-371.

[MFK] D. Mumford, J. Fogarty and F. Kirwan, Geometric invariant theory, Springer Verlag, 1994.

[PV] M. van der Put and H. Voskuil, Symmetric spaces associated to split algebraic groups over a local field, J. f. d. reine u. angew. Math. 433 (1992),

69 - 100.

[Ra] M. Rapoport, Period domains over finite and local fields, Proc. of Symp. in Pure Math. 62.1 (1997), 361-381.

[RZ] M. Rapoport and Th. Zink, Period spaces for p-divisible groups, Princeton University Press, 1996.

[Ro] M. Ronan, Lectures on buildings, Academic Press Inc., 1989.

[S] T.A. Springer, Linear algebraic groups, Birkhäuser, 1981.

[T.1] J. Tits, Classification of algebraic semisimple groups, Proc. of symp. in pure math. 9 (1966), 33-62.

[T.2] J. Tits, Représentations linéaires irréductibles d'un groupe réductif sur un corps quelconque, J. f. d. reine u. angew. Math. 247 (1971), 196-220.

[T.3] J. Tits, Reductive groups over local fields, Proc. A.M.S. Symp. Pure Math. 33 (1979), 29-69.

[To] B. Totaro, Tensor products in p-adic Hodge theory, Duke Math. J. 83 (1996), 79-104.

[V.1] H. Voskuil, Nonarchimedean flag domains, preprint M.P.I. 92-84 (1992).

[V.2] H. Voskuil, Nonarchimedean flag domains and semistability, preprint university of Sydney 1998-35,

(http://www.maths.usyd.edu.au:8000/res/Analysis/tr.html)